\documentclass[10pt]{amsart}
\usepackage{amssymb,amsmath,amsfonts,amscd,color,euscript,graphics,comment}
\usepackage{enumerate}
\usepackage{amssymb}
\usepackage{mathtools}
\usepackage{amsthm}
\usepackage{latexsym}
\usepackage{multicol}
\usepackage{verbatim}
\usepackage[all]{xy}
\usepackage{tikz}
\usepackage{hyperref}

\advance\textwidth by 1.2in \advance\oddsidemargin by -.6in \advance\evensidemargin by -.6in

\newtheorem{cor}{Corollary}[subsection]
\newtheorem{lem}{Lemma}[subsection]
\newtheorem{prop}{Proposition}[subsection]

\theoremstyle{definition}
\newtheorem*{defn}{Definition}

\theoremstyle{definition}
\newtheorem{thm}{Theorem}[subsection]
\newtheorem*{theom}{Theorem}
\newtheorem*{conje}{Conjecture}

\newtheorem{rem}{Remark}[subsection]

\newcounter{cnt}
 \makeatletter
\def\mydggeometry{\makeatletter\dg@YGRID=1\dg@XGRID=20\unitlength=0.003pt\makeatother}
\makeatother \theoremstyle{remark}

\numberwithin{equation}{section}

\let\bwdg\bigwedge
\def\bigwedge{{\textstyle\bwdg}}
\begin{document}
\newcommand{\al}{\alpha}
\newcommand{\alvee}{{\al^\vee}}
\newcommand{\bevee}{{\be^\vee}}
\newcommand{\phivee}{{\Phi^\vee}}
\newcommand{\Om}{{\Omega}}
\newcommand{\g}{\gamma}
\newcommand{\G}{\Gamma}
\newcommand{\BSlam}{\hat\Sigma(\g_\lam )}
\newcommand{\de}{\delta}
\newcommand{\ta}{\theta}
\newcommand{\De}{\Delta}
\newcommand{\Supp}{\operatorname{Supp}}
\newcommand{\semi}{{\,\rule[.1pt]{.4pt}{5.3pt}\hskip-1.9pt\times}}
\newcommand{\M}{{\mathcal M}}
\newcommand{\A}{\mathbb A}
\newcommand{\bfb}{\mathbf b}
\newcommand{\charc}{\hbox{\rm Char}\,}
\newcommand{\vep}{\varepsilon}
\newcommand{\set}[1]{\left\{#1\right\}}
\newcommand{\si}{\sigma}
\newcommand{\avee}{\alpha^\vee}
\newcommand{\bvee}{\beta^\vee}
\newcommand{\Xvee}{{X^\vee}}
\newcommand{\Cg}{{\mathfrak{C}(\Lhg)}}
\newcommand{\Xveep}{{X^\vee_+}}
\newcommand{\Phivee}{{\Phi^\vee}}
\newcommand{\Phiveep}{{\Phi^\vee_+}}
\newcommand{\Phiveeplus}{{\Phi^\vee_+}}
\newcommand{\Rvee}{{R^\vee}}
\newcommand{\Rveep}{{R^\vee_+}}
\newcommand{\Delvee}{\Delta^\vee}
\newcommand{\Waff}{{W^{{\mathfrak a}}}}
\newcommand{\Saff}{{S^{{\mathfrak a}}}}
\newcommand{\Aaff}{{A^{{\mathfrak{a}}}}}
\newcommand{\eps}{\epsilon}
\newcommand{\om}{\omega}
\newcommand{\wt}{\widetilde}
\newcommand{\delp}{\partial^+}
\newcommand{\delm}{\partial^-}
\newcommand{\lam}{{\lambda}}
\newcommand{\chp}{{\check{P}}}
\newcommand{\Lam}{{\Lambda}}
\newcommand{\wte}{{\widetilde{e}}}
\newcommand{\Lhg}{\widehat{\mathfrak{g}}}
\newcommand{\Lhh}{\widehat{\mathfrak{h}}}
\newcommand{\Lhb}{\widehat{\mathfrak{b}}}
\newcommand{\Lhn}{\widehat{\mathfrak{n}}}
\newcommand{\wDelta}{\widehat{{\Delta}}}
\newcommand{\wPhi}{\widehat{{\Phi}}}
\newcommand{\ph}{\widehat{P}}
\newcommand{\waff}{{W^{\rm aff}}}
\newcommand{\wiff}{{\widetilde{W}^{\rm aff}}}
\newcommand{\QRightarrow}{%
\begingroup
\tikzset{every path/.style={}}%
\tikz \draw (0,3pt) -- ++(1em,0) (0,1pt) -- ++(1em+1pt,0) (0,-1pt) -- ++(1em+1pt,0) (0,-3pt) -- ++(1em,0) (1em-1pt,5pt) to[out=-75,in=135] (1em+2pt,0) to[out=-135,in=75] (1em-1pt,-5pt);
\endgroup
}

%%%%%%%%%%%%%%%%%%%%%%
%%%%%Gothic%%%%%%%%%%%%%%%%%
%%%%%%%%%%%%%%%%%%%%%%
\newcommand{\La}{{\mathfrak{a}}}
\newcommand{\Lgg}{{\mathfrak{g}}}
\newcommand{\Lb}{{\mathfrak{b}}}
\newcommand{\Lh}{{\mathfrak{h}}}
\newcommand{\Ln}{{\mathfrak{n}}}
\newcommand{\Ld}{{\mathfrak{d}}}
\newcommand{\LU}{{\mathfrak{U}}}
\newcommand{\Lu}{{\mathfrak{u}}}
\newcommand{\ug}{{\underline{G}}}
\newcommand{\Lc}{{\mathfrak{C}}}
\newcommand{\Ls}{{\mathfrak{s}}}
\newcommand{\Lcinf}{{\mathfrak{C}_{-\infty}}}
\newcommand{\Lgc}{{\mathfrak{g} \otimes \mathbb{C}[t]}}
\newcommand{\Lgl}{{\mathfrak{g} \otimes \mathbb{C}[t,t^{-1}]}}
\newcommand{\Lbc}{{\mathfrak{I_b}}}
\newcommand{\thmref}[1]{Theorem~\ref{#1}}
\newcommand{\secref}[1]{Section~\ref{#1}}
\newcommand{\ssecref}[1]{~\ref{#1}}
\newcommand{\lemref}[1]{Lemma~\ref{#1}}
\newcommand{\propref}[1]{Proposition~\ref{#1}}
\newcommand{\corref}[1]{Corollary~\ref{#1}}
\newcommand{\remref}[1]{Remark~\ref{#1}}
\newcommand{\defref}[1]{Definition~\ref{#1}}
\newcommand{\er}[1]{(\ref{#1})}
\newcommand{\id}{\operatorname{id}}
\newcommand{\Lg}{\lie g[t]}
\newcommand{\Lggs}{\lie g^{\Gamma}}
\newcommand{\Lgs}{\lie g[t]^{\sigma}}
\newcommand{\Lgn}{\mathfrak{n^{-}}[t]^{\sigma}}
\newcommand{\Lgb}{\mathfrak{b}[t]^{\sigma}}
\newcommand{\Nsl}{\mathfrak{n^{-}_{sl_2}}}
\newcommand{\Ws}{W^{\Gamma}(n\omega)}
\newcommand{\ord}{\operatorname{\emph{ord}}}
\newcommand{\tensor}{\otimes}
\newcommand{\from}{\leftarrow}
\newcommand{\nc}{\newcommand}
\newcommand{\rnc}{\renewcommand}
\newcommand{\dist}{\operatorname{dist}}
\newcommand{\qbinom}[2]{\genfrac[]{0pt}0{#1}{#2}}
\newcommand{\C}{\ensuremath{\mathbb{C}}}
\newcommand{\N}{\ensuremath{\mathbb{N}}}
\newcommand{\Z}{\ensuremath{\mathbb{Z}}}
\newcommand{\Q}{\ensuremath{\mathbb{Q}}}
\nc{\cal}{\mathcal} \nc{\goth}{\mathfrak} \rnc{\bold}{\mathbf}
\renewcommand{\frak}{\mathfrak}
\renewcommand{\Bbb}{\mathbb}
\newcommand{\lie}[1]{\mathfrak{#1}}
 \nc{\Hom}{\operatorname{Hom}} \nc{\op}{\operatorname}
  \nc{\mode}{\operatorname{mod}}
 \nc{\Mod}{\operatorname{Mod}}  
\nc{\End}{\operatorname{End}} \nc{\wh}[1]{\widehat{#1}} \nc{\Ext}{\operatorname{Ext}} \nc{\ch}{\text{ch}} \nc{\ev}{\operatorname{ev}}
\nc{\Ob}{\operatorname{Ob}} \nc{\soc}{\operatorname{soc}} \nc{\rad}{\operatorname{rad}} \nc{\head}{\operatorname{head}}
\nc{\orb}{\operatorname{orb}}
\def\und{\underline}
\nc{\ci}{{\mathcal I}}

 \nc{\Cal}{\cal} \nc{\Xp}[1]{X^+(#1)} \nc{\Xm}[1]{X^-(#1)}
\nc{\on}{\operatorname}\nc{\J}{{\cal J}} \renewcommand{\P}{{\cal
P}} \nc\boa{\bold a} \nc\bob{\bold b} \nc\boc{\bold c} \nc\bod{\bold d} \nc\boe{\bold e} \nc\bof{\bold f}
\nc\bog{\bold g} \nc\boh{\bold h} \nc\boi{\bold i} \nc\boj{\bold j} \nc\bok{\bold k} \nc\bol{\bold l} \nc\bom{\bold m}
\nc\bon{\bold n} \nc\boo{\bold o} \nc\bop{\bold p} \nc\boq{\bold q} \nc\bor{\bold r} \nc\bos{\bold s} \nc\boT{\bold t}
\nc\boF{\bold F} \nc\bou{\bold u} \nc\bov{\bold v} \nc\bow{\bold w} \nc\boz{\bold z} \nc\boy{\bold y} \nc\ba{\bold A}
\nc\bb{\bold B} \nc\bc{\bold C} \nc\bd{\bold D} \nc\be{\bold E} \nc\bg{\bold G} \nc\bh{\bold H} \nc\bi{\bold I} \nc\bj{\bold
J} \nc\bk{\bold K} \nc\bl{\bold L} \nc\bm{\bold M} \nc\bo{\bold O} \nc\bp{\bold P} \nc\bq{\bold Q}
\nc\br{\bold R} \nc\bs{\bold S} \nc\bt{\bold T} \nc\bu{\bold U} \nc\bv{\bold V} \nc\bw{\bold W} \nc\bx{\bold
x} \nc\KR{\bold{KR}} \nc\rk{\bold{rk}} \nc\het{\text{ht }}

\nc\toa{\tilde a} \nc\tob{\tilde b} \nc\toc{\tilde c} \nc\tod{\tilde d} \nc\toe{\tilde e} \nc\tof{\tilde f} \nc\tog{\tilde g}
\nc\toh{\tilde h} \nc\toi{\tilde i} \nc\toj{\tilde j} \nc\tok{\tilde k} \nc\tol{\tilde l} \nc\tom{\tilde m} \nc\ton{\tilde n}
\nc\too{\tilde o} \nc\toq{\tilde q} \nc\tor{\tilde r} \nc\tos{\tilde s} \nc\toT{\tilde t} \nc\tou{\tilde u} \nc\tov{\tilde v}
\nc\tow{\tilde w} \nc\toz{\tilde z}

\title
{Demazure modules and Weyl modules: The twisted current case}
\author{Ghislain Fourier and Deniz Kus}
\address{Ghislain Fourier:\newline
Mathematisches Institut, Universit\"at zu K\"oln, Germany}
\email{gfourier@math.uni-koeln.de}
\address{Deniz Kus:\newline
Mathematisches Institut, Universit\"at zu K\"oln, Germany}
\email{dkus@math.uni-koeln.de}
\begin{abstract}
We study finite--dimensional respresentations of twisted current algebras and show that any graded twisted Weyl module is isomorphic to level one Demazure module for the twisted affine Kac-Moody algebra. Using the tensor product property of Demazure modules, we obtain, by analyzing the fundamental Weyl modules, dimension and character formulas. Moreover we prove that graded twisted Weyl modules can be obtained by taking the associated graded modules of Weyl modules for the loop algebra, which implies that its dimension and classical character are independent of the support and depend only on its classical highest weight. These results were known before for untwisted current algebras and are new for all twisted types.
\end{abstract}

\maketitle \thispagestyle{empty}
%%%%%%%%%%%%%%%%%%%%%%%%%%%%%%%%%%%%%%%%%%%%%%%%%%%%%%%%%%%%%%%%%%%%%%%%%%%%%%%%%%%%%%%%%%%%%%%%%%%%%%%%%%%%%%%%%%%%%%%%%%%%%%%%%%%
%         Introduction
%%%%%%%%%%%%%%%%%%%%%%%%%%%%%%%%%%%%%%%%%%%%%%%%%%%%%%%%%%%%%%%%%%%%%%%%%%%%%%%%%%%%%%%%%%%%%%%%%%%%%%%%%%%%%%%%%%%%%%%%%%%%%%%%%%%
\section{Introduction}
Weyl modules for loop algebras $\lie g \otimes \C[t,t^{-1}]$, where $\Lgg$ is a simple complex Lie algebra, have gained a lot of attraction during the last two decades. Starting with the analysis of finite--dimensional irreducible modules for quantum affine algebras (\cite{CP01}), which are highest weight modules in a certain sense. It was natural to ask for maximal finite--dimensional modules with these highest weights since contrary to the theory of simple complex Lie algebras, the category of finite--dimensional modules is not semi--simple. 
In the same paper it was conjectured, that the classical limit $q = 1$ of these irreducible modules specialize to modules for the loop algebra satisfying some universal properties, the so called \textit{local Weyl modules}. In a series of papers (\cite{BN04}, \cite{CL06}, \cite{CP2001}, \cite{FoL07}, \cite{Nak01}, \cite{Na11}) the character and dimension of these Weyl modules were computed. 
In the proofs, these modules were identified with Weyl modules for the current algebra $\lie g \otimes \C[t]$. Using the tensor product property (\cite{CP01}) and some pullback maps, the study was reduced to analyzing graded Weyl modules for $\lie g \otimes \C[t]$, where the grading is induced by the grading of $\C[t]$.\\
One major step in the analysis of the graded Weyl modules is their identification with level one Demazure modules for simply--laced algebras (\cite{CL06}, \cite{FoL07}). With the tensor product property for Demazure modules (\cite{FoL06}) and the computation for fundamental Weyl and Demazure modules (\cite{CL06}, \cite{FoL06}), the character and dimension formulas were proven. 
In the non simply--laced case, Weyl modules admit a filtration by Demazure modules and via this filtration, the dimension and character formula were proven (\cite{Na11}). One should mention that these results can also be deduced from the results in \cite{BN04}, \cite{Nak01}, but there is no written proof so far in the literature.\\
Local Weyl modules for current and loop algebras can be parametrized by finitely supported functions from $\C$ (resp. $\C^*$) to $P^+$, the set of dominant integral weights for $\Lgg$. To each function one can associate a weight, which is the sum of all images, hence in $P^+$. 
To summarize the results above, the dimension and character of a local Weyl module are independent of the support of the parametrizing function and depend only on its weight. The graded local Weyl module of weight $\lambda$ is parametrized by the function of weight $\lambda$ with support in the origin only. We can also reformulate this result in terms of the global Weyl module, which is a projective module in a certain category and in general infinite--dimensional. The results on local Weyl modules are equivalent to the statement, that the global Weyl module is a free module for a certain commutative algebra $\mathbf{A}_{\lambda}$.\\
There are several ways to generalize the notion of local Weyl modules. By replacing $\C[t,t^{-1}]$ with a commutative, associative algebra (\cite{FL04}, \cite{CFK10}) one can define local and global Weyl modules as before, obtain similar tensor product properties, but character and dimension formulas are known only in certain cases. Even for a case as simply looking as $\lie g = \lie{sl}_2$ and $\C[t_1,\ldots, t_n]$ with $n \geq 4$ there is no dimension formula known.\\
Another way of generalizing local Weyl modules is to look at twisted current and loop algebras. Given a complex simple Lie algebra $\lie g$ and a commutative algebra $A$ ($=\C[t], \C[t,t^{-1}]$), both equipped with the action of a finite group $\Gamma$ $(\Gamma=\Z/m\Z)$ by automorphism, one can extend this action to $\lie g \otimes A$. The fixpoint Lie algebra $(\lie g \otimes A)^{\Gamma}$ is called the twisted current algebra (resp. twisted loop algebra). 
The twisted current algebra is a subalgebra of the twisted affine Kac-Moody algebra associated to $\lie g$, while the twisted loop is obtained by taking the quotient by the central element of the subalgebra without derivation \cite{Ca05}.\\
Local Weyl modules for the twisted loop algebra were introduced and studied in \cite{CFS08}. It was proven, that every Weyl module is the tensor product of Weyl modules located in a single point only. So to obtain dimension and character formulas it was sufficient to compute them for Weyl modules with support in a single point. 
The main theorem in \cite{CFS08} states that every Weyl module for the twisted loop algebra is isomorphic to the restriction of a Weyl module for the untwisted loop algebra. So all interesting information can be deduced from this isomorphism. In \cite{FMS} the aforementioned global Weyl modules will be defined and studied for twisted loop algebras as well. It will be shown, that the twisted global Weyl module is a submodule of the untwisted global Weyl module, viewed as a module for the twisted loop algebra by restriction. The results about twisted local Weyl module translate again into the freeness of the twisted global Weyl module as a module for a certain commutative algebra $\mathbf{A}_{\lambda}^{\Gamma}$.\\
In \cite{FKKS11} the notion of local Weyl modules was generalized to certain equivariant map algebras. 
Given $X$ an affine scheme and $\lie g$ a finite--dimensional Lie algebra, both defined over an algebraically closed field and $\Gamma$ a finite group acting on $X$ and $\lie g$ by automorphisms, the equivariant map algebra is the Lie algebra of equivariant maps from $X$ to $\lie g$. In \cite{FKKS11} several restrictions to this general case were assumed, the group action on $X$ had to be free and abelian. But under these assumptions, again the tensor product property was proven. Furthermore it was shown, that every Weyl module for the equivariant map algebra is isomorphic to the restriction of a Weyl module for the algebra of maps from $X$ to $\lie g$.\\
In this paper we are considering the gap in the computation of dimension and character formulas for local Weyl modules of twisted current algebras. For twisted and untwisted loop and current algebras, dimension formulas for all local Weyl modules are known except for graded local Weyl modules for the twisted current algebra. Let $\Gamma$ be the finite group of non--trivial diagram automorphism of a simple Lie algebra $\lie g$, so $\Gamma$ is of order $2$ or $3$ and $\lie g$ of type $A, D, E$. In terms of equivariant map algebras, the affine scheme would be $X = \C$ and $\Gamma=\langle\xi\rangle$, where $\xi$ is the multiplication by a primitive $2$nd or $3$rd root of unity. We see immediately that $0$ is a fix point, so the group action is not free. In this setting, the results of \cite{FKKS11} do not apply at the origin. \\
The goal of this paper is to compute a dimension and character formula for the local Weyl module located in $0$ (the graded local Weyl module) of the twisted current algebra. The main tool are, as in \cite{FoL07} and \cite{Na11}, Demazure modules.\\
There are two cases to be considered, the first one is: 
\begin{theom}
Let $\lie g$ be not of type $A_{2l}$ and $\lambda \in P_0^+$, then the local graded $(\Lgg\otimes\C[t])^{\Gamma}$-Weyl module $W^{\Gamma}(\lambda)$ is isomorphic to a Demazure module of level $1$.
\end{theom}
In the proof we will use the $\lie{sl}_2 \otimes \C[t]$ and the $(\lie{sl}_3 \otimes \C[t])^\Gamma$ cases (proven in \cite{CP01}, \cite{CL06}, resp. Section~\ref{section7}). A tensor product property for Demazure modules was proven in \cite{FoL06}, so to obtain a character formula for Weyl modules it is sufficient to determine the \textit{fundamental local Weyl modules}, as done in Section~\ref{section5}. Concluding we were able to prove an analogous result to \cite{CFS08}, \cite{FKKS11}, that the dimension of the local Weyl module does not depend on the support but only on the highest weight.
\begin{theom} 
For $\lie g$ not of type $A_{2l}$ and $\lambda \in P_0^+$, the local graded $(\Lgg\otimes\C[t])^{\Gamma}$-Weyl module is isomorphic to the associated graded module of the restriction of a local Weyl module for $\Lgg\otimes\C[t,t^{-1}]$.
\end{theom}
In the second case, we assume that $\lie g$ is of type $A_{2l}$, then the fixpoint algebra $\lie g_0$ is of type $B_l$. Here with our methods, one can only determine the local Weyl module for weights $\lambda$, where $\lambda(\alpha_l^{\vee})$ is odd. 
In this case there is an identification with Demazure modules as before, so the graded local $(\Lgg\otimes\C[t])^{\Gamma}$-Weyl module is isomorphic to a Demazure module of level $1$. Furthermore we are able to show the following: 
\begin{theom}
Let $\lambda=\lambda_1+\lambda_2\in P^+_0$, where $\lambda_2(\alpha_l^{\vee})$ is odd, and $a\in\C^*$. Then $$W^{\Gamma}(\lambda)\cong \hbox{\rm gr}(W_{a}(\lambda_1)\otimes W^{\Gamma}(\lambda_2)),$$ where $W_{a}(\lambda_1)$ is the local Weyl module for $\Lgg\otimes\C[t,t^{-1}]$, supported in $a$ with highest weight $\lambda_1$.
\end{theom}

In the case where $\lambda(\alpha_l^{\vee})$ is even the dimension and character of the local Weyl modules remains uncomputed, the identification with Demazure modules fails. We can state here a conjecture only
\begin{conje}
Let $\lambda \in P_0^+$, then the graded local Weyl module is isomorphic to the associated graded module of the restriction of a local Weyl module for $\Lgg\otimes\C[t,t^{-1}]$. The dimension of a local Weyl module of highest weight $\lambda$ is independent of the support of the module. 
\end{conje}
The structure of the paper is as follows, in Section~\ref{section2} are basics and notations for affine Kac-Moody algebras recalled, in Section~\ref{section3} for twisted current algebras. In Section~\ref{section4} Demazure and Weyl modules are defined. In Section~\ref{section5} we identify Demazure modules with Weyl modules and determine the ``smallest'' Weyl modules. 
In Section~\ref{section6} we show that every graded Weyl module of the twisted current algebra can be obtained by taking the associated graded of the restriction of a untwisted loop module. In Section~\ref{section7}, the case $\lie g = \lie{sl}_3$ is treated seperately, since it is used in some of the proofs of the other cases.
\vskip12pt
\textbf{Acknowledgements:} 
We would like to thank Vyjayanthi Chari and Peter Littelmann for helpful discussions. We would also like to thank the Hausdorff Research Institute for Mathematics and the organizers of the Trimester Program on the Interaction of Representation Theory with Geometry and Combinatorics, during which the ideas in the current paper were developed. The first author was partially sponsored by the DFG-Schwerpunktprogramm 1388 ``Darstellungstheorie" and the second author by the ``SFB/TR 12-Symmetries and Universality in Mesoscopic Systems".
%%%%%%%%%%%%%%%%%%%%%%%%%%%%%%%%%%%%%%%%%%%%%%%%%%%%%%%%%%%%%%%%%%%%%%%%%%%%%%%%%%%%%%%%%%%%%%%%%%%%%%%%%%%%%%%%%%%%%%%%%%%%%%%%%%%
%         The affine Kac-Moody algebras
%%%%%%%%%%%%%%%%%%%%%%%%%%%%%%%%%%%%%%%%%%%%%%%%%%%%%%%%%%%%%%%%%%%%%%%%%%%%%%%%%%%%%%%%%%%%%%%%%%%%%%%%%%%%%%%%%%%%%%%%%%%%%%%%%%%
\section{The affine Kac-Moody algebras}\label{section2}
\subsection{Notation and basic results}
In this section we fix the notation and the usual technical padding.
Let $\Lgg=\Lgg(A)$ be a simple complex Lie algebra of rank $l$ associated to a Cartan matrix $A$ of finite type. We fix a Cartan subalgebra
$\Lh$ in $\Lgg$ and a Borel subalgebra $\Lb\supseteq \Lh$. Denote
$\Phi\subseteq \Lh^*$ the root system of $\Lgg$, and, corresponding to the choice
of $\Lb$, let $\Phi^+$ be the set of positive roots
and let $\Pi=\{\al_1,\ldots,\al_l\}$ be the corresponding basis of $\Phi$.\vskip4pt
For a root $\beta\in\Phi$ let $\beta^\vee\in\Lh$ be its coroot.
The basis of the dual root system (also called the coroot system)
$\Phi^\vee\subset \Lh$ is denoted $\Pi^\vee=\{\al_1^\vee,\ldots,\al_l^\vee\}$.
The Weyl group $W$ of $\Phi$ is generated by the simple reflections $s_i=s_{\al_i}$ associated to the simple roots.\vskip4pt
Let $P=\bigoplus^{l}_{i=1}\Z \omega_i$ be the weight lattice of $\Lgg$ and let $P^+=\bigoplus^{l}_{i=1}\Z_{\geq 0} \omega_i$ be the subset
of dominant weights.  The group algebra of $P$ is denoted $\Z[P]$,
we write $\chi=\sum a_\mu e^\mu$ (finite sum, $\mu\in P$, $a_\mu\in\Z$)
for an element in $\Z[P]$, where the embedding $P\hookrightarrow \Z[P]$ is
defined by $\mu\mapsto e^\mu$. Further we denote by $Q=\bigoplus^{l}_{i=1}\Z \al_i$ 
(respectively $Q^{+}=\bigoplus^{l}_{i=1}\Z_{\geq} \al_i$)
be the root (respectively positive root) lattice and let $\{x^{\pm}_i,h_i|i\in I\}$ be a set of Chevalley generators of $\Lgg$.\vskip4pt
Let $\Lhg$ be the affine Kac--Moody algebra (twisted or untwisted) corresponding
to the Cartan matrix $\widehat{A}=(a_{i,j})$. Note that, if $\Lhg$ is a untwisted affine Kac--Moody algebra associated to $\Lgg$:
$$
\Lhg=\Lgg\otimes_\C\C[t,t^{-1}]\oplus \C c\oplus \C d.
$$
Here $d$ denotes the derivation $d=t\frac{d}{dt}$ and $c$ is
the canonical central element. Recall that the Lie bracket is
defined as
$$
[t^m\otimes x + \lam c +\mu d, t^n\otimes y + \nu c + \eta d]=
t^{m+n}\otimes [x,y] + \mu n t^n\otimes y +\eta m t^m\otimes x +
m\delta_{m,-n}(x, y) c.
$$
We assume $\Lhg$ is arbitrary (possibly twisted) and we fix a Cartan subalgebra $\Lhh$ in $\Lhg$ and a Borel subalgebra $\Lhb\supseteq \Lhh$, $\Pi=\{\al_0,\ldots,\al_l\}$ the set of simple roots, $\Pi^\vee=\{\al_0^\vee,\ldots,\al_l^\vee\}$ the set of simple coroots. Denote by $\widehat\Phi$ the root system of $\Lhg$ and
let $\wPhi^+$ be the subset of positive roots. We denote by $\widehat{P}$ the weight lattice of $\Lhg$ and let $\widehat{P^+}$ be the subset of dominant weights. The Weyl group $\widehat{W}$ of $\widehat{\Phi}$ is generated
by the simple reflections $s_i=s_{\al_i}$ associated to the simple roots.
Further we fix uniquely determined vectors $w=(a_0,\ldots,a_l)^t, v=(a_0^\vee,\dots,a_l^\vee)$, such that $v\widehat{A}=\widehat{A}w=0.$ Then it is known that the center of $\Lhg$ is 1-dimensional and is spanned by the \emph{canonical cantral element} $$c=\sum_{i=0}^{l}a_i^\vee \al_i^\vee.$$ Define further $$\delta=\sum_{i=0}^{l}a_i \al_i;\quad \theta=\delta-a_0\al_0$$ and $d\in\Lhh$ which satisifes the following conditions $$\al_i(d)=0, \mbox{for } i=1,\ldots,l ;\quad \al_0(d)=1.$$ Clearly the elements $\al_0^\vee,\ldots,\al_l^\vee,d$ form a basis of $\Lhh$.
We have a non--degenerate symmetric bilinear
form $\langle\cdot,\cdot\rangle$ on $\Lhh$ defined in (\cite{K90}, Chapter 6)
\begin{equation}\label{hatform}
\left\{
\begin{array}{ll}
\langle\al_i^\vee,\al_j^\vee\rangle=\frac{a_j}{a_j^\vee}a_{i,j}&i,j=0,\ldots,\ell\\
\langle\al_i^\vee,d\rangle=0&i=1,\ldots,\ell\\
\langle\al_0^\vee,d\rangle=\frac{a_0}{a_0^\vee}&\langle d,d\rangle=0.\\
\end{array}\right.
\end{equation}
This $\widehat{W}$-invariant form induces a isomorphism
$$
\nu:\Lhh \longrightarrow \Lhh^*, \quad
\nu(h):\left\{
\begin{array}{rcl}
\Lhh&\rightarrow&\C\\
h'&\mapsto &\langle h,h'\rangle\\
\end{array}\right.
$$
With the notation as above it follows
for $i=0,\ldots,l$:
$$
\nu(\al_i^\vee)=\frac{a_i}{a_i^\vee}\al_i
$$

Let $\Lam_0,\ldots,\Lam_l$ be the fundamental weights in $\widehat{P^+}$,
then for $i=1,\ldots,l$ we have
\begin{equation}\label{Lamdecomp}
\Lam_i=\om_i+\frac{a_i^\vee}{a_0^{\vee}}\Lam_0.
\end{equation}
With this we have $\widehat{P}=\sum_{i=0}^{l}\Z \Lambda_i + \Z(\delta/a_0)$ and $\widehat{P^+}=\sum_{i=0}^{l}\Z_{\geq 0}\Lambda_i + \Z(\delta/a_0)$.

\subsection{Realisation of twisted affine algebras}
In this paper we are mainly interested in twisted affine Kac-Moody algebras, which can be realised as fixed point subalgebras of so-called twisted graph automorphisms. Let $\Lgg$ be a finite dimensional simple Lie algebra and $\sigma:\Lgg\rightarrow \Lgg$ be a graph automorphism of order $m$. In particular $$m =\left\{\begin{array}{cl} 2, & \mbox{if $\Lgg$ of type $A_{2l},A_{2l-1}, D_{l+1}$ or $E_6$ }  \\ 3, & \mbox{if $\Lgg$ is of type $D_4$}\end{array} \right. $$ Let $\xi$ be a primitive $m^{\mbox{th}}$ root of unity, then it is well-known that there exists a decomposition of $\Lgg$ into eigenspaces. We obtain: 
$$\Lgg=\Lgg_{0}\oplus\cdots\oplus \Lgg_{m-1},$$
whereby $\Lgg_{j}=\{x\in \Lgg | \sigma(x)= \xi^j x\}$, $j=0,\cdots ,m-1$. The fixed point algebra $\Lgg_{0}$ is again a simple complex Lie algebra of type $C_l,B_l,F_4$ or $G_2$ and the eigenpaces are irreducible $\Lgg_{0}$-modules. 
\begin{rem}
Let $\mathfrak{a}$ be a subalgebra of $\Lgg$ such that $\sigma(\mathfrak{a})=\mathfrak{a}$, then we get a analogue decomposition $$\mathfrak{a}=\mathfrak{a}_{0}\oplus\cdots \oplus \mathfrak{a}_{m-1}.$$
So if $\Lgg=\Ln\oplus\Lh\oplus\Ln_{-}$ is a triangular decomposition of $\Lgg$, we obtain $$\Lgg_{j}=\Ln_{j}\oplus\Lh_{j}\oplus(\Ln_{-})_{j}\mbox{ for all $0\leq j \leq m-1$}.$$
\end{rem}
Now we can extend $\sigma$ to a automorphism of the corresponding untwisted affine algebra given by $$\sigma(x\otimes t^{i})=\xi^{-i} \sigma(x)\otimes t^{i} \quad \mbox{for }x\in\Lgg$$
$$\sigma(c)=c;\quad \sigma(d)=d.$$
The twisted affine algebra is realized as the fixed point subalgebra $$\Lhg\cong\bigoplus_{k\in\Z}(\Lgg_{0}\otimes t^{mk})\oplus\cdots\oplus \bigoplus_{k\in\Z}(\Lgg_{m-1}\otimes t^{mk+(m-1)})\oplus \C c\oplus \C d$$
$$=\bigoplus_{j=0}^{m-1}\bigoplus_{k\in\Z}(\Lgg_j \otimes t^{mk+j})\oplus \C c \oplus \C d.$$ 
Using the above notation we can conclude $$\Lhh=\Lh_{0}\oplus(\C c + \C d)\quad \Lhh^*=(\Lh_{0})^*\oplus(\C \delta + \C \Lambda_{0}).$$ 
We have the following table, which describes the various possibilities for $\Lgg, \Lgg_{0}$, $\Lhg$ and the eigenspaces $\Lgg_{1}, \Lgg_{2}.$
$$\begin{tabular}[ht]{|l|l|c|c|c|c|c|}
  \hline
   $m$ & $\Lgg$ & $\lie g_{0}$ & $\Lhg$& $\Lgg_{1}$& $\Lgg_{2}$& Dynkin diagram of $\Lhg$\\
  \hline\hline
      2 & $A_{2}$ & $A_1$ & $A^{(2)}_{2}$& $V(4\omega_{1})$& $/$&$\underset{\mathclap{0}}{\circ} \QRightarrow \underset{\mathclap{1}}{\circ}$\\
    2 & $A_{2l}, l \geq 2$ & $B_l$ & $A^{(2)}_{2l}$& $V(2\omega_{1})$& $/$&$\underset{\mathclap{0}}{\circ} \Rightarrow \underset{\mathclap{1}}{\circ} - \dotsb - \underset{\mathclap{l-1}}{\circ} \Rightarrow \underset{\mathclap{l}}{\circ}$\\
 2 & $A_{2l-1}, l \geq 2$ & $C_l$ & $A^{(2)}_{2l-1}$& $V(\omega_{2})$& $/$&$\underset{\mathclap{1}}{\circ} - \underset{\mathclap{2}}{\overset{\overset{\textstyle\circ_{\mathrlap{0}}}{\textstyle\vert}}{\circ}} \,-\, \dotsb - \underset{\mathclap{l-1}}{\circ}\Leftarrow \underset{\mathclap{l}}{\circ}$\\
  2 & $D_{l+1}, l\geq 3$ & $ B_l$ & $D^{(2)}_{l+1}$& $V(\omega_{1})$& $/$&$\underset{\mathclap{0}}{\circ} \Leftarrow \underset{\mathclap{1}}{\circ} - \dotsb - \underset{\mathclap{l-1}}{\circ} \Rightarrow \underset{\mathclap{l}}{\circ}$\\
  2 & $E_6$ & $F_4$ & $E^{(2)}_{6}$& $V(\omega_{1})$& $/$&$\underset{\mathclap{0}}{\circ}-\underset{\mathclap{1}}{\circ} - \underset{\mathclap{2}}{\circ} \Leftarrow \underset{\mathclap{3}}{\circ} - \underset{\mathclap{4}}{\circ}$\\
  3 & $D_4$& $G_2$ & $D^{(3)}_{4}$& $V(\omega_{2})$& $V(\omega_{2})$&$ \underset{\mathclap{1}}{\circ} \Rrightarrow \underset{\mathclap{2}}{\circ}-\underset{\mathclap{0}}{\circ}$\\
  \hline
\end{tabular}$$
We put a ``0" on (almost) everything related to $\Lgg_0$, e.g. denote by $\Phi_0\subseteq (\Lh_0)^*$ the root system of $\Lgg_0$.
The recently defined element $\delta$ is the imaginary root in $\wPhi^+$ and $\theta$ is the highest short root of the root system of $\Lgg_{0}$ if $\widehat{A}$ is of type $A^{(2)}_{2l-1}, D^{(2)}_{l+1}, E^{(2)}_{6}, D^{(3)}_{4}.$ In the remaining twisted cases $\theta-\al_1$ is the highest root of the root system of $\Lgg_0$.
For more details we refer to (\cite{K90},\cite{Ca05}).

\begin{rem}
The untwisted Kac-Moody algebras $\Lhg=\Lgg\otimes_\C\C[t,t^{-1}]\oplus \C c\oplus \C d$ can also be realised as fixed point algebras for any automorphism of order 1. We have $\Lgg_0=\Lgg$ and the eigenspaces are the zerospaces. In this case $\theta$ is the highest root of $\Lgg$.
\end{rem}

\subsection{The extended affine Weyl group}
Now we give a description of the Weyl group $\widehat{W}$ of the affine Kac-Moody algebra $\Lhg$. The Weyl group is generated by fundamental reflections $s_0,\ldots,s_l$, which act on $\Lhh^*$ by $$s_i(\lambda)=\lambda-\lambda(\al_i^\vee)\al_i, \quad \lambda\in\Lhh^*.$$ Since $\delta(\al_i^\vee)=0$ for all $i=0,\ldots,l$, the Weyl group $\widehat{W}$ fixes $\delta.$ Another well-known description of the affine Weyl-group is the following. Let $W_0$ be the subgroup of $\widehat{W}$  generated by $s_1,\ldots,s_l$, i.e. $W_0$ can be identified with the Weyl group of the Lie algebra $\Lgg_0$, since $W_0$ operates trivially on $(\C \delta + \C \Lambda_{0})$. Further let \begin{equation}\label{weylgrouptrans1}M=\sum_{i=1}^{l} \Z \al_i\quad \mbox{if }\widehat{A} \mbox{ symmetric or } m > a_0 \end{equation} or \begin{equation}\label{weylgrouptrans2}  M=\nu(\sum_{i=1}^{l} \Z \al^\vee_i)\quad \mbox{otherwise}.\end{equation} For an element $\mu\in M$ let $t_\mu$ the following endomorphism of the vector space $\Lhh^{*}$: 
%%%%%%%%%%%%%%%%%%%%%%%%%%%%%%%%%%%%%%%%%%
\begin{equation}\label{muaction}
\Lam =\lam +b\Lam_0+r\delta\mapsto
t_\mu(\Lam)=\Lam+\Lam(c)\mu-(\langle\Lam,\mu\rangle+\frac{1}{2}\langle\mu,\mu\rangle\Lam(c))\delta
\end{equation}
Obviously we have $t_{\mu}\circ t_{\mu'}=t_{\mu+\mu'}$, denote $t_M$ the
abelian group consisting of the elements
$t_\mu$, $\mu\in M$. Then $\widehat{W}$  is the semidirect product
$\widehat{W} = W_0 \semi t_M$.

The {\it extended affine Weyl group} $\widehat{W}^{ext}$ is the semidirect product
$\widehat{W}^{ext}= \widehat{W}\semi t_L$, where $L=\nu(\bigoplus_{i=1}^l\Z\om^\vee_i)$
is the image of the coweight lattice. The action of an element $t_\mu$,
$\mu\in L$, is defined as above in (\ref{muaction}).
Let $\Sigma$ be the subgroup of $\widehat{W}^{ext}$ stabilizing the dominant Weyl chamber
$\widehat{C}$:
$$
\Sigma=\{\sigma\in\widehat{W}^{ext}\mid \sigma(\widehat{C})=\widehat{C}\}.
$$
Then $\Sigma$ provides a complete system of coset representatives of
$\widehat{W}^{ext}/\widehat{W}$, so we can write in fact $\widehat{W}^{ext}=\Sigma\semi\widehat{W}$.

The elements $\sigma\in\Sigma$ are all of the form
$$
\sigma=\tau_i t_{-\nu(\omega_i^\vee)}= \tau _i t_{-\omega_i},
$$
where $\omega_i^\vee$ is a minuscule fundamental coweight. Further,
set $\tau_i=w_0 w_{0,i}$, where $w_0$ is the longest word $W_0$ and
$w_{0,i}$ is the longest word in $W_{\om_i}$, the stabilizer of $\om_i$ in $W_0$.

\subsection{Weight space decomposition and roots}
Remember that the Borel subalgebra for the twisted case is given by:$$\Lhb=((\Lh_0\oplus\Ln_0)\otimes1)\oplus\bigoplus_{k\in\N_{>0}}(\Lgg_0\otimes t^{mk})\oplus\bigoplus_{k\in\N}(\Lgg_{1}\otimes t^{mk+1})\oplus\cdots\oplus \bigoplus_{k\in\N}(\Lgg_{m-1}\otimes t^{mk+(m-1)})\oplus \C c\oplus \C d.$$  
Furthermore we remember that $\Lgg_{j}$ is a irreducible $\Lgg_0$-module for all $j$, so one can obtain the following weight space decomposition 
$$\Lgg_{j}=\bigoplus_{\alpha\in(\Lh_{0})^*}(\Lgg_{j})_{\alpha}$$

\begin{prop}
$\Lh_{j}=(\Lgg_{j})_0;\quad (\Ln_{-})_{j}= \bigoplus_{\alpha\in \Phi^{-}\vline_{\Lh_{0}}}(\Lgg_{j})_{\alpha};\quad \Ln_{j}= \bigoplus_{\alpha\in \Phi^+\vline_{\Lh_0}}(\Lgg_{j})_{\alpha}$, $0 \leq j \leq m-1$.
\end{prop}
Let $(\Phi_0)_s$ be the set of short roots and  $(\Phi_0)_l$ be the set of long roots of  $\lie g_0$ and $\Phi_j=\{\al\in(\mathfrak{h}_0)^{*}|(\lie g_j)_{\al}\neq 0\}-\{0\}$, then we get $$\Lgg_{j}=\Lh_j\oplus\bigoplus_{\alpha\in \Phi_j}(\Lgg_{j})_{\alpha}=\Lh_j\oplus\bigoplus_{\alpha\in \Phi^+_j}(\Lgg_{j})_{\alpha}\oplus\bigoplus_{\alpha\in \Phi^-_j}(\Lgg_{j})_{\alpha},$$ whereby $\dim(\lie g_j)_{\al}=1$ for all $\al\in \Phi_j$, hence $(\lie g_j)_{\pm\al}=\C X^{\pm}_{\alpha,j}$ for $\al\in \Phi^{+}_j$ and we have the following table \cite{Ca05}:
$$\begin{tabular}[ht]{|l|l|c|c|c|}
  \hline
    $\mathfrak{g}$ & $\Lgg_0$ & $\Phi_1$ & $\Phi_2$& Dynkin diagram of $\Lgg_0$\\
  \hline\hline
   $A_{2}$ & $A_1$ & $(\Phi_0)\cup\{2\alpha : \alpha\in\Phi_0\}$& $/$&$\underset{\mathclap{1}}{\circ}$ \\
   $A_{2l}, l\geq 2$ & $B_l$ & $(\Phi_0)\cup\{2\alpha : \alpha\in(\Phi_0)_s\}$& $/$&$\underset{\mathclap{1}}{\circ} - \underset{\mathclap{2}}{\circ} - \dotsb - \underset{\mathclap{l-1}}{\circ} \Rightarrow \underset{\mathclap{l}}{\circ}$ \\
   $A_{2l-1}, l \geq 2$ & $C_l$ & $(\Phi_0)_s$ & $/$&$\underset{\mathclap{1}}{\circ} - \underset{\mathclap{2}}{\circ} - \dotsb - \underset{\mathclap{l-1}}{\circ} \Leftarrow \underset{\mathclap{l}}{\circ}$ \\
   $D_{l+1}, l\geq 3$ & $B_l$ & $ (\Phi_0)_s$ & $/$&$\underset{\mathclap{1}}{\circ} - \underset{\mathclap{2}}{\circ} - \dotsb - \underset{\mathclap{l-1}}{\circ} \Rightarrow \underset{\mathclap{l}}{\circ}$\\
   $E_6 $ & $F_4$ & $ (\Phi_0)_s$ & $/$&$\underset{\mathclap{1}}{\circ} - \underset{\mathclap{2}}{\circ} \Leftarrow \underset{\mathclap{3}}{\circ} - \underset{\mathclap{4}}{\circ}$\\
   $D_{4}$ & $G_2$ & $ (\Phi_0)_s$ & $(\Phi_0)_s$&$ \underset{\mathclap{1}}{\circ} \Rrightarrow \underset{\mathclap{2}}{\circ}$\\
   
  \hline
\end{tabular}$$
%%%%%%%%%%%%%%%%%%%%%%%%%%%%%%%%%%%%%%%%%%%%%%%%%%%%%%%%%%%%%%%%%%%%%%%%%%%%%%%%%%%%%%%%%%%%%%%%%%%%%%%%%%%%%%%%%%%%%%%%%%%%%%%%%%%
%         
%%%%%%%%%%%%%%%%%%%%%%%%%%%%%%%%%%%%%%%%%%%%%%%%%%%%%%%%%%%%%%%%%%%%%%%%%%%%%%%%%%%%%%%%%%%%%%%%%%%%%%%%%%%%%%%%%%%%%%%%%%%%%%%%%%%
\section{The twisted current algebra \texorpdfstring{$\Cg$}{}}\label{section3}
In this section we will define the twisted current algebra $\Cg$ and certain subalgebras, which will be needed in the following sections. The main object of this paper will be $$\Cg:=\bigoplus_{j=0}^{m-1}\bigoplus_{k\geq 0}(\Lgg_j \otimes t^{mk+j}).$$
The algebra $\Cg$ can be realized by taking the fixpoints under the group of automorphisms $\Gamma$ restricted to the current algebra, in detail $(\Lgg\otimes \C[t])^{\Gamma}\cong \Cg$, hence it is called the \emph{twisted current algebra}.\\

In order to give an explicit basis of $\Cg$ we use the embedding $\Lgg_j\hookrightarrow \Lgg$\ for all $0\leq j\leq m-1$, so that we can realize the generators of the weight spaces $(\Lgg_j)_{\pm\alpha}$ as elements in $\Lgg$. This is already described in \cite{K90},\cite{Ca05} and \cite{CFS08} if $\alpha$ is a simple root and can be continued to arbitrary $\alpha\in\Phi_0$: Let $(\widetilde{\alpha_1},\cdots,\widetilde{\alpha_m})$ be a m-element orbit of $\sigma$ on $\Phi$ and $x^{\pm}_{\widetilde{\alpha_i}}\in \Lgg$ be root vectors such that $\sigma(x^{\pm}_{\widetilde{\alpha_i}})=x^{\pm}_{\widetilde{\alpha_{i+1}}}.$ Then we obtain  $$(\sum^{m-1}_{i=0}(\xi^i)^jx^{\pm}_{\sigma^i(\widetilde{\alpha_1})})\in(\Lgg_j)_{\pm\widetilde{\alpha_1}\vline_{\Lh_{0}}},\ 0\leq j\leq m-1.$$
In (\cite{Ca05}, Chapter 18.4) it is shown that the weight spaces of $\Lgg_j$ are spanned by such elements for all m-element orbits $(\widetilde{\alpha_1},\cdots,\widetilde{\alpha_m})$. So the weight spaces $(g_j)_{\pm\alpha}$ can be described as follows: There has to be a root $\widetilde{\alpha}$, such that $\widetilde{\alpha}\vline_{\Lh_{0}}=\alpha$ and \begin{equation}\label{basisro}\C X^{\pm}_{\alpha,j}= \C(\sum^{m-1}_{i=0}(\xi^i)^jx^{\pm}_{\sigma^i(\widetilde{\alpha})})\end{equation}
We set further 
\begin{equation}\label{basisco}\C h_{\alpha,j}=\C(\sum^{m-1}_{i=0}(\xi^i)^jh_{\sigma^i(\widetilde{\alpha})})\end{equation}
At this point we have adapted our notation while we denote by $h_{\alpha,0}$ the coroot of a root $\alpha\in \Phi_0$.
\begin{lem}\label{isosl}
Assume $\Lhg$ is of type $A^{(2)}_{2l-1}, D^{(2)}_{l+1}, E^{(2)}_{6}$ or $D^{(3)}_{4}$. If $\al$ is a long root then we get an canonical isomorphism $$\lie{sl}_2\otimes \C[t]\cong\langle X^{\pm}_{\al,0}\otimes t^{ms}, h_{\al,0}\otimes t^{ms}| s\in\N\rangle_{\C} =: \lie{sl}_{2,\al}\otimes \C[t^m]$$ and if $\al$ is short we have $$\lie{sl}_2\otimes \C[t]\cong\langle X^{\pm}_{\al,j}\otimes t^{ms+j}, h_{\al,j}\otimes t^{ms+j}| s\in\N\ ,\ 0\leq j \leq m-1\rangle_{\C}=:\lie{sl}_{2,\al}\otimes \C[t].$$
\proof
Since the Lie algebra $\langle X^{\pm}_{\al,0}, h_{\al,0}\rangle_{\C}$ is canonically isomorph to $\lie{sl}_2$ the first isomorphism is given by $$x^{\pm}\otimes t^s \mapsto X^{\pm}_{\al,0}\otimes t^{ms}$$$$h\otimes t^s \mapsto h_{\al,0}\otimes t^{ms}.$$ To verify the second isomorphism we define $$x^{\pm}\otimes t^s \mapsto X^{\pm}_{\al,j}\otimes t^{s},\ \mbox{if $s\equiv j\mod m$}$$$$h\otimes t^s \mapsto h_{\al,j}\otimes t^{s},\ \mbox{if $s\equiv j\mod m$}$$ 
To show that this map is an homomorphism of Lie algebras we need to check 
\begin{equation}\label{checklie}[X^{+}_{\al,i_1},X^{-}_{\al,i_2}]= h_{\al,i_1+i_2\mod m},\ [h_{\al,i_2},X^{\pm}_{\al,i_1}]=\pm 2 X^{\pm}_{\al,i_1+i_2\mod m}\end{equation}
Since we require $\alpha$ to be a short root, we know that the weight space $(g_j)_{\pm\alpha},\ 0\leq j \leq m-1$ is non-zero and therefore we can use the description in (\ref{basisro}), (\ref{basisco}) with $\widetilde{\alpha}$, such that $\sigma(\widetilde{\alpha})\neq\widetilde{\alpha}$. More than this, a case by case consideration shows $\sigma^j(\widetilde{\alpha})(\sigma^i(\widetilde{\alpha})^{\vee})=0$ and $\sigma^j(\widetilde{\alpha})-\sigma^i(\widetilde{\alpha})$ is not a root for $i\neq j$, e.g. in type $D^{(2)}_{l+1}$ we have for an arbitrary short root $\al_i+\cdots+\al_l$ of $B_l$, that $\widetilde{\al}=\al_i+\cdots+\al_l$ and therefore $\sigma(\widetilde{\al})-\widetilde{\al}=\al_{l+1}-\al_l$ is not a root and $\sigma(\widetilde{\al})(\widetilde{\al}^{\vee})=\widetilde{\al}(\sigma(\widetilde{\al})^{\vee})=0$. The proof in the other cases is similar. We set $X^{\pm}_{\alpha,j}=(\sum^{m-1}_{i=0}(\xi^i)^jx^{\pm}_{\sigma^i(\widetilde{\alpha})}),\ h_{\alpha,j}=(\sum^{m-1}_{i=0}(\xi^i)^j h_{\sigma^i(\widetilde{\alpha})})$. The required equations in (\ref{checklie}) are now immediate.
\endproof
\end{lem}
If $\Lhg$ is of type $A^{(2)}_{2l}$ we obtain a similar result
\begin{lem}\label{isoa2}
Assume $\Lhg$ is of type $A^{(2)}_{2l}$ and $\al$ be a long root then we get an canonical isomorphism $$\lie{sl}_2\otimes \C[t]\cong\langle X^{\pm}_{\al,j}\otimes t^{ms+j}, h_{\al,j}\otimes t^{ms+j}| s\in\N\ ,\ 0\leq j \leq m-1\rangle_{\C}=:\lie{sl}_{2,\alpha}\otimes\C[t]$$ and if $\al$ is a short root, then we get an canonical isomorphism  $$\mathfrak{C}(A^{(2)}_{2}) \cong \langle X^{\pm}_{\al,j}\otimes t^{ms+j}, X^{\pm}_{2\al,1}\otimes t^{ms+1}, h_{\al,j}\otimes t^{ms+j}|s\in \N ,\ 0\leq j \leq m-1\rangle_{\C}.$$
\proof
The proof of the first isomorphism is similar to Lemma~\ref{isosl} and to justify the second isomorphism we will demonstrate how to realize the elements $h_{\al,j}, X^{\pm}_{\al,j}$, $X^{\pm}_{2\al,1}$ as elements in $A_{2l}$. Let $\al=\al_i+\cdots+\al_l$ be an arbitrary short root of type $B_l$ and $\widetilde{\al}=\al_i+\cdots+\al_l$ be the root considered as a root in type $A_{2l}$, i.e. the restriction to $\Lh_0$ equals $\al$. It is easy to see that $\sigma(\widetilde{\al})\neq\widetilde{\al}$, $\sigma(\widetilde{\al})-\widetilde{\al}$ is not a root of $A_{2l}$ and continuing $\sigma(\widetilde{\al})(\widetilde{\al}^{\vee})=\widetilde{\al}(\sigma(\widetilde{\al})^{\vee})=-1$. We set $$X^{\pm}_{\al,j}=(\xi)^j\sqrt{2}(x^{\pm}_{\widetilde{\al}}+\xi^jx^{\pm}_{\sigma(\widetilde{\al})})\in (\Lgg_j)_{\pm\al}$$
$$X^{\pm}_{2\al,1}=[x^{\pm}_{\widetilde{\al}},x^{\pm}_{\sigma(\widetilde{\al})}]\in (\Lgg_1)_{\pm 2\al}$$
$$h_{\al,j}=2^{\delta_{0,j}}(h_{\widetilde{\al}}+\xi^jh_{\sigma(\widetilde{\al})})$$
Now, knowing the embedding in $A_{2l}$, it is straighforward to check the required relations.
\endproof
\end{lem}
\subsection{Filtration on $\Cg$}
The Lie algebra $\Cg$ has a natural grading and an associated natural filtration $F^{\bullet}(\Cg)$, where $F^{s}(\Cg)$ is defined
to be the subspace of $\Lgg$-valued polynomials with degree smaller or equal $s$.
One has an induced filtration also on the enveloping algebra $U(\Cg)$ and therefore an induced filtration on arbitrary cyclic $U(\Cg)$-modules $W$ with cyclic vector $w$.
Denote by
$W_{s}$ the subspace spanned by the vectors of the form $g.w$,
where $g \in \, F^{s}(U(\Cg))$, and denote the associated graded $\Cg$--module
by $\hbox{\rm gr}(W)$
$$
\hbox{\rm gr}(W) = \bigoplus_{i\ge 0} W_i / W_{i-1},\ \mbox{ where }W_{-1} = 0.
$$

%%%%%%%%%%%%%%%%%%%%%%%%%%%%%%%%%%%%%%%%%%%%%%%%%%%%%%%%%%%%%%%%%%%%%%%%%%%%%%%%%%%%%%%%%%%%%%%%%%%%%%%%%%%%%%%%%%%%%%%%%%%%%%%%%%%
%         
%%%%%%%%%%%%%%%%%%%%%%%%%%%%%%%%%%%%%%%%%%%%%%%%%%%%%%%%%%%%%%%%%%%%%%%%%%%%%%%%%%%%%%%%%%%%%%%%%%%%%%%%%%%%%%%%%%%%%%%%%%%%%%%%%%%
\section{Demazure modules and Weyl modules}\label{section4}
\subsection{Definition of Demazure modules}
For a dominant weight $\Lam\in \ph^+$ let $V(\Lam)$
be the irreducible highest weight module of highest weight $\Lam$. Given an element
$w\in\widehat{W}$, fix a generator $v_{w(\Lam)}$
of the line $V(\Lam)_{w(\Lam)}=\C  v_{w(\Lam)}$
of $\Lhh$--eigenvectors in $V(\Lam)$ of weight $w(\Lam)$.
\begin{defn}
The $U(\Lhb)$--submodule $V_w(\Lam)=U(\Lhb).v_{w(\Lam)}$
generated by $ v_{w(\Lam)}$ is called the {\it Demazure submodule
of $V(\Lam)$} associated to $w$.
\end{defn}
\begin{rem}\mbox{}
\begin{enumerate}
\item Since $\Lhh$ acts by multiplication with a scalar on $v_{w(\Lam)}$, the Demazure module $V_w(\Lam)$ is a cyclic
$U(\Lhn)$--module generated.
\item The modules $V_w(\Lam)$ are finite--dimensional although $V(\Lambda)$ is infinite--dimensional.
\end{enumerate}
\end{rem}

To associate more generally to every element $\sigma w\in\widehat{W}^{ext}=\Sigma\semi\widehat{W}$
a Demazure module, recall that elements in $\Sigma$ correspond to automorphisms
of the Dynkin diagram of $\Lhg$, and thus define an associated automorphism
of $\Lhg$, also denoted $\sigma$. For a module $V$ of $\Lhg$
let $\sigma^{*}(V)$ be the module with the twisted action $g\circ v=\sigma^{-1}(g)v$.
Then for the irreducible module of highest weight $\Lam\in \ph^+$ we get
$\sigma^{*}(V(\Lam))=V(\sigma(\Lam))$.

So for $\sigma w\in \widehat{W}^{ext}=\Sigma\semi\widehat{W}$ we set
\begin{equation}\label{demazuredefn}
V_{\sigma w}(\Lam):=V_{\sigma w \sigma^{-1}}(\sigma(\Lam)).
\end{equation}

We are mainly interested in $\Lgg_0$-stable Demazure modules. For $i\in I_0$ we have $X^{-}_{\al_i,0}v_{w(\Lam)}=0$ if and only if $w(\Lam)(\al_i^\vee)\leq 0$. Consequently we can see that $V_{w}(\Lam)$ is $\Lgg_0$-stable if and only if $w(\Lam)(\al_i^\vee)\leq 0$ for all $i\in I_0$.
Assume that $\omega(\Lambda)=-\lambda+k\Lambda_0+i\delta$, then $V_{w}(\Lam)$ is stable under $\Lgg_0$ if and only if $\lambda\in P_0^+$. 
We define a set $$X=\{(\lambda,k,i)\in P_0^+\times(1/a^{\vee}_0)\Z_{>0}\times(1/a_0)\Z \ |\ \exists! \Lambda\in \widehat{P}^+: w_0(\lambda)+k\Lambda_0+i\delta\in\widehat{W}(\Lambda) \},$$ where $w_0$ is the longest word in $W_0$. Let $(\lambda,k,i)\in X$ and $\omega\in\widehat{W}$, such that $\omega(\Lambda)=w_0(\lambda)+k\Lambda_0+i\delta$. Then by the above computation we get the $\Lgg_0$-stability of the Demazure module $V_{w}(\Lam)$ and we denote $$V_{w}(\Lam)=D(k,\lambda)[i].$$
\begin{rem}\mbox{}
\begin{enumerate}
\item The $\Lgg_0$ stable Demazure modules are in fact $\Cg$-modules.
\item For any $\Lam\in \ph^+$ and $i\in(1/a_0)\Z$, we have $V(\Lambda)\cong V(\Lambda+i\delta)$, as $\Cg$-modules. Therefore we get $$D(k,\lambda)[i]\cong D(k,\lambda)[i+n],$$ which justifies the notation $D(k,\lambda)$ as a $\Cg$-module.
\end{enumerate}
\end{rem}
\begin{rem}
Whenever we speak about  $D(k,\lambda)$ we will assume that $(\lambda,k)\in X$. If $\Lhg$ is not of type $A^{(2)}_{2l}\ (l\geq 1)$ the set $X$ is given by $X=P_0^+\times \Z_{>0}\times \Z$ and else we have $P_0^+\times \Z_{>0}\times \Z \subsetneq X$.
\end{rem}

\subsection{Demazure character formula}
Let $\beta$ be a real root of the root system $\widehat\Phi$. We define
the {\it Demazure operator:}
$$
D_{\beta}:\Z[\ph]\rightarrow\Z[\ph],\quad
D_{\beta}(e^\lam)=\frac{e^\lam-e^{s_\beta(\lam)-\beta}}{1-e^{-\beta}}
$$
\begin{lem}\label{demazurezerleg}\mbox{}
\begin{enumerate}
\item For $\lambda,\mu\in\ph$ we have:
\begin{equation}\label{DemazureExpl}
D_{\beta}(e^\lam)=\left\{
\begin{array}{ll}
e^\lam+e^{\lam-\beta}+\dots+e^{s_\beta(\lam)}&
\hbox{\rm if\ }\langle\lam,\beta^\vee\rangle\ge 0\\
0&\hbox{\rm if\ }\langle\lam,\beta^\vee\rangle=-1\\
-e^{\lam+\beta}-e^{\lam+2\beta}-\dots-e^{s_\beta(\lam)-\beta}&
\hbox{\rm if\ }\langle\lam,\beta^\vee\rangle\le -2\\
\end{array}
\right.
\end{equation}
\item Let $\chi,\eta\in\Z[\ph]$. If $D_{\beta}(\eta)=\eta$, then
\begin{equation}\label{chimalcharformel}D_{\beta}(\chi\cdot\eta)=\eta\cdot(D_{\beta}(\chi)).\end{equation}
\end{enumerate}
\end{lem}
\proof For (1) see (\cite{D74}, (1.5)--(1.8)) and for (2) see (\cite{FoL06}, (2.2)).
\endproof
\vskip 5pt
Since $D_{\alpha_i}(1-e^\delta)=(1-e^\delta)$ for all $i=0,\ldots,n$,
(\ref{chimalcharformel}) shows that the ideal $I_\delta=\langle (1-e^\delta)\rangle$ is
stable under all Demazure operators $D_{\beta}$. Thus we obtain
induced operators (we still use the same notation $D_\beta$)
$$
D_{\beta}:\Z[\ph]/I_\delta\longrightarrow \Z[\ph]/I_\delta,\quad
e^\lam + I_\delta\mapsto D_{\beta}(e^\lam) + I_\delta.
$$
In the following we denote by $D_i$, $i=0,\ldots,n$ the Demazure operator
$D_{\al_i}$ corresponding to the simple root $\al_i$. Recall that for any reduced decomposition $w=s_{i_1}\cdots s_{i_r}$
of $w\in\widehat{W}$ the operator $D_w=D_{i_1}\cdots D_{i_r}$
is independent of the choice of the decomposition (see \cite{KU02},
Corollary 8.2.10). We have the following important theorem:
\begin{thm}[\cite{KU02} Chapter VIII]\label{demazurcharacterformula}
$$
\charc V_w(\Lam)=D_{w}(e^\Lam).
$$
\end{thm}
We will need the following elementary proposition:
\begin{prop}\label{opaus}
Let $\lambda_1^{\vee}, \lambda_2^{\vee}$ be two dominant coweights, and set $\lambda^{\vee}=\lambda_1^{\vee}+\lambda_2^{\vee}$. Then \begin{enumerate}
\item $D_{t_{-\nu(\lambda_1^{\vee})}}D_{t_{-\nu(\lambda_2^{\vee})}}=D_{t_{-\nu(\lambda^{\vee})}}$
\item $D_{t_{-\nu(\lambda_1^{\vee})}}D_{\omega_0}=D_{t_{-\nu(\lambda_1^{\vee})}\omega_0}$
\end{enumerate}
\end{prop}

\subsection{Properties of Demazure modules}
Since $V_w(\Lam)=U(\Lhb)\cdot v_{w(\Lam)}$, there exists an Ideal $J\subseteq U(\Lhb)$, such that $V_w(\Lam)\cong U(\Lhb)/J$. So the Demazure module can be described by generators and relations, which was done in \cite{M88}. We give here a reformulation for the twisted affine case:

\begin{prop}[\cite{M88}]
Let $\Lam\in \ph^+$ and let $w$ be an element of the affine
Weyl group of $\Lhg$. The Demazure module $V_w(\Lam)$ is as a $U(\Lhb)$-module
isomorphic to the cyclic module, generated by $v \neq 0$ with respect to the following
relations.\\ 
For $\beta\in\Phi^{+}_j$, $0\leq j \leq m-1$ we have:
\begin{eqnarray*}
&&(X^{+}_{\beta,j} \otimes t^{ms+j})^{k_{\beta}+1}.v = 0 \quad \mbox{ where } s\ge 0,\quad k_{\beta} =
max\{0,-\langle w(\Lam), (\beta+(ms+j)\delta)^{\vee}\rangle\} \\
&&(X_{\beta,j}^{-} \otimes t^{ms+j})^{k_{\beta}+1}.v = 0 \quad \mbox{ where } s>-\delta_{j,\{1,\cdots,m-1\}},\quad k_{\beta} =
max\{0,-\langle w(\Lam),(-\beta + (ms+j)\delta)^{\vee}\rangle\} \\
&&(h\otimes t^{ms+j}).v =\delta_{j,0}\delta_{s,0} w(\Lam)(h) v\quad \forall h\in \Lh_j,\quad \mbox{ where } s\ge 0,\quad d.v=w(\Lam)(d).v,\ c.v=w(\Lam)(c) v
\end{eqnarray*}
\end{prop}

\begin{cor}\label{demazurerelations}
As a module for $\Cg$ the Demazure module
$D(k, \lam)$ is isomorphic to the cyclic $U(\Cg)$--module generated by a
vector $v\neq 0$ subject to the following relations:\\
For $\beta\in\Phi^{+}_j$, $0\leq j \leq m-1$ we have:
\begin{eqnarray*}
&&\Ln_{j} \otimes t^{j}\C[t^{m}].v = 0\\
&&(X^{-}_{\beta,j} \otimes t^{ms+j})^{k_{\beta}+1}.v = 0\quad \mbox{where } s\ge 0,\quad k_{\beta} =
\max\{ 0 ,\langle\lambda, \beta^{\vee}\rangle-\frac{2(ms+j)}{\langle\beta,\beta\rangle}ka_{0}^{\vee} \}\\
&&(h\otimes t^{ms+j}).v =\delta_{j,0}\delta_{s,0} \lam(h) v\quad \forall h\in \Lh_j,\quad \mbox{where } s\ge 0\end{eqnarray*}
\proof
The proof is similar to the one given in (\cite{FoL07} Corollary 1).
\endproof
\end{cor}

\begin{rem} Since the defining relations of $D(k, \lam)$ respect the grading of $\Cg$, $D(k, \lam)$ is a graded module.
\end{rem}
\vskip5pt

In \cite{FoL06} it was shown by using the Demazure operator, that $D(k, \lam)$ decomposes as a $\lie g$ (resp. $\lie g_0$) module into a tensor product of ''smaller'' Demazure modules. We give here the result for the twisted affine case:

\begin{thm}\cite{FoL06}\label{tensorzerlegung} Let $\lambda^{\vee} = \lambda^{\vee}_1 + \lambda^{\vee}_2 +
\ldots + \lambda^{\vee}_r$ be a sum of dominant coweights. Then for $m \geq 0$ we have an
isomorphism of $\Lgg_0$-modules between the Demazure module
$V_{-\lambda^{\vee}}(m\Lam_0)$ and the tensor product of
Demazure modules:
$$
V_{{-\lambda^\vee}}(m\Lam_0)\simeq
V_{{-\lambda_1^\vee}}(m\Lam_0)\otimes
V_{{-\lambda_2^\vee}}(m\Lam_0) \otimes\cdots\otimes
V_{{-\lambda_r^\vee}}(m\Lam_0).
$$
\end{thm}
\begin{rem} This theorem holds for any special vertex $k$ of the twisted affine diagram. 
\end{rem}%%%%%%%%%%%%%%%%%%%%%%%%%%%%%%%%%%%%%%%%%%%%%%%%%%%%%%%%%%%%%%%%%%%%%%%%%%%%%%%%%%%%%%%%%%%%%%%%%%%%%%%%%%%%%%%%%%%%%%%%%%%%%%%%%%%
%         
%%%%%%%%%%%%%%%%%%%%%%%%%%%%%%%%%%%%%%%%%%%%%%%%%%%%%%%%%%%%%%%%%%%%%%%%%%%%%%%%%%%%%%%%%%%%%%%%%%%%%%%%%%%%%%%%%%%%%%%%%%%%%%%%%%%
\subsection{Definition of Weyl modules}
The representation theory of twisted current algebras is particularly interesting because the category of finite--dimensional representation is not semisimple. It makes sense to ask for the ``maximal" finite--dimensional cyclic representations in this class, which leads to the definition of Weyl modules.
Let $\lambda=\sum^{l}_{i=1} m_i\omega_i\in P^{+}_{0}$ be a dominant integral weight for $\Lgg_0$. Then we define the Weyl module $W^{\Gamma}(\lambda)$ in terms of generators and relations:
\begin{defn}\label{weyldef}
Let $\lambda=\sum^{l}_{i=1} m_i\omega_i$ be a dominant integral weight for $\Lgg_0$. Denote $W^{\Gamma}(\lambda)$ be the $U(\Cg)$-module generated by an element $w_{\lambda}$ with the relations:
\begin{equation}\label{weylrel} \Ln_j\otimes t^j\C[t^m].w_{\lambda}=0,\ 0\leq j\leq m-1\end{equation}
\begin{equation}\label{weylrel2}(h\otimes t^{ms+j}).w_{\lambda} =\delta_{j,0}\delta_{s,0} \lam(h) w_{\lambda}\quad \forall h\in \Lh_j,\quad  \mbox{where } s\ge 0\end{equation}
\begin{equation}\label{weylrel3}(X^{-}_{\beta,0}\otimes 1)^{\lambda(\beta^{\vee})+1}w_{\lambda}=0,\ \mbox{for all positive roots $\beta$ of $\Lgg_0$}\end{equation}
\end{defn}
\begin{rem}
Note that the modules  $W^{\Gamma}(\lambda)$ are graded modules since $U(\Cg)$ is graded by the powers of $t$ and the defining relations are graded, particulary we have $$W^{\Gamma}(\lambda)\cong\bigoplus_{s\in\Z_{+}}W^{\Gamma}(\lambda)[s],$$ where $W^{\Gamma}(\lambda)[s]$ is a $\Lgg_0$-module by identifying $\Lgg_0$ with $\Lgg_0\otimes 1\subseteq \Cg$.
\end{rem}
\subsection{Properties of Weyl modules}
\begin{prop}\label{weylprop}\mbox{}
\begin{enumerate}
\item We have $$W^{\Gamma}(\lambda)=\bigoplus_{\mu\in(\Lh_0)^*}W^{\Gamma}(\lambda)_\mu$$ and $W^{\Gamma}(\lambda)_\mu\neq 0$ only if $\mu\in \lambda-Q^{+}_0$. Further we get $W^{\Gamma}(\lambda)_\mu\neq 0$ if and only if $W^{\Gamma}(\lambda)_{w(\mu)}\neq 0$ for all $w\in W_0.$
\item As a $\Lgg_0$ module $W^{\Gamma}(\lambda)$ and $W^{\Gamma}(\lambda)[s]$ decompose into finite--dimensional irreducible representations of $\Lgg_0$.
\item Let $\mu$ be a dominant integral weight, such that $\lambda-\mu$ is as well dominant integral. Then there exists a canonical homomorphism $W^{\Gamma}(\lambda)\rightarrow W^{\Gamma}(\mu)\otimes W^{\Gamma}(\lambda-\mu)$ mapping $w_{\lambda}$ to $w_{\mu}\otimes w_{\lambda-\mu}.$
\end{enumerate}
\proof
It sufficies to show that for every $v\in W^{\Gamma}({\lambda})_\mu$ the module $U(\lie g_0).v$ is finite dimensional, since this proves the non-trivial statements in part (1) and (2). Part (3) is clear from the defining relations. Given $v\in W^{\Gamma}({\lambda})_\mu$ we obtain $\bold U(\lie g_0).v=\bold U((\Ln_{-})_0)\bold U(\Ln_0).v.$ From part (1) we obtain that $U(\Ln_0).v$ is finite dimensional. By the PBW-theorem $U((\Ln_{-})_0)$ is spanned by monomials, so it suffices to show that $X^{-}_{\beta,0}\in(\Ln_{-})_0$ acts  nilpotently on v. Assume that $v\in U(\Cg)w_{\lambda}$ and 
the action of $(\Ln_{-})_0$ on $\Cg$, which is given by the Lie bracket is locally nilpotent.
We obtain with $$(X^{-}_{\beta,0}\otimes 1)^{\lambda(\beta^\vee)+1}w_{\lambda}=0,\quad (X^{-}_{\beta,0})^N(\underbrace{u.w_{\lambda}}_{=v})=\sum_{k=0}^{N}\binom{N}{k}((X^{-}_{\beta,0})^{k}u) (X^{-}_{\beta,0})^{N-k} w_{\lambda}$$ that $X^{-}_{\beta,0}$ acts nilpotently on v, which finally implies that $\bold U(\lie g_0).v$ is finite dimensional. 
\endproof
\end{prop}
\begin{rem}
$W^{\Gamma}(\lambda)$ is finite--dimensional. This will be an immediate consequence of Theorem~\ref{mainthm} and Corollar~\ref{fid}.
\end{rem}

By definition we obtain some obvious maps between Weyl modules and certain Demazure modules.
\begin{cor}\label{weyldemquot}
Let $\lambda$ be a dominant integral weight for $\Lgg_0$. Then for all $k\in (1/a_0^{\vee})\Z_{>0}$, such that $(\lambda,k)\in X$, the Demazure module $D(k,\lambda)$ is a quotient of the Weyl module $W^{\Gamma}(\lambda)$.
\proof
This follows immediately by comparing the relations for the Weyl module in Definition~\ref{weyldef} and the relations for the Demazure module in Corollary~\ref{demazurerelations}.
\endproof
\end{cor}

In this paper we want to show, that the map between Weyl  and Demazure modules is in fact an isomorphism. This is already known for untwisted current algebras of simply-laced type (\cite{CP01}, \cite{CL06},\cite{FoL07}). We recall the result for $\lie g = \lie{sl}_2$ here only, since this will be heavily used throughout this paper.
\begin{thm}\label{sl2isobek}
For $\lie g = \lie{sl}_2$ and $n \omega \in P^+$, we have an isomorphism of $\lie{sl}_2\otimes\C[t]$-modules
$$W(n \omega) \cong D(1, n\omega).$$
\end{thm}

%%%%%%%%%%%%%%%%%%%%%%%%%%%%%%%%%%%%%%%%%%%%%%%%%%%%%%%%%%%%%
%      Connection between Weyl modules and Demazure modules
%%%%%%%%%%%%%%%%%%%%%%%%%%%%%%%%%%%%%%%%%%%%%%%%%%%%%%%%%%%%%

\section{Connection between Weyl modules and Demazure modules}\label{section5}
In this section we will show, that almost all Weyl modules are isomorphic to certain Demazure modules, e.g. the map in Corollary~\ref{weyldemquot} is in fact an isomorphism.
\begin{thm}\label{mainthm}
Suppose $\Lhg$ is of type $A^{(2)}_{2l-1}, D^{(2)}_{l+1}, E^{(2)}_{6}$ or $D^{(3)}_{4}$, then we have an isomorphism of $\Cg$-modules 
$$W^{\Gamma}(\lambda)\cong D(1/a^{\vee}_0,\lambda).$$ 
If $\Lhg$ is of type $A^{(2)}_{2l}$ and $\lambda = \sum\limits_{i=1}^{l} m_i \omega_i$ be a dominant weight, such that $m_l$ is odd, we have an isomorphism of $\Cg$-modules 
$$W^{\Gamma}(\lambda)\cong D(1/a^{\vee}_0,\lambda).$$
\proof
By Corollary~\ref{weyldemquot} we know already that the Demazure module is a quotient of the Weyl module. By comparing the defining relations in Corollary~\ref{demazurerelations} and in Definition~\ref{weyldef}, we see that to prove that this map is
an isomorphism, it is sufficient to show that the generator of the Weyl module is subject to the following relations:\\
For all $0 \leq j \leq m-1$, $\beta\in\Phi^{+}_j$:
\begin{equation}\label{pot1}(X^{-}_{\beta,j} \otimes t^{ms+j})^{k_{\beta}+1}.w_{\lambda} = 0,\ \mbox{where } s\ge 0,\ k_{\beta} =
\max\{ 0 ,\langle\lambda, \beta^{\vee}\rangle-\frac{2(ms+j)}{\langle\beta,\beta\rangle}\frac{1}{a_{0}^{\vee}}a_{0}^{\vee} \}.\end{equation}
Assume $\Lhg$ is not of type $A^{(2)}_{2l}$, then (\ref{pot1}) is equivalent to :
\begin{equation}\label{pot2}(X^{-}_{\beta,j} \otimes t^{ms+j})^{k_{\beta}+1}.w_{\lambda} = 0,\ \mbox{where } s\ge 0,\ k_{\beta} =
\left\{\begin{array}{cl} \max\{ 0 ,\langle\lambda, \beta^{\vee}\rangle-s \}, & \mbox{if $\beta$ is long}\\ \max\{ 0 ,\langle\lambda, \beta^{\vee}\rangle-(ms+j) \}, & \mbox{if $\beta$ is short} \end{array}\right.\end{equation}
Let $\beta\in \Phi^+_0$ be a long root and $V=U(\lie{sl}_{2,\beta}\otimes \C[t^m]).w_{\lambda}\subseteq W^{\Gamma}(\lambda)$ be the $\lie{sl}_{2,\beta}\otimes \C[t^m]$-submodule. Further let $W(\langle\lambda,\beta^{\vee}\rangle\omega)$ be the $\lie{sl}_2\otimes \C[t]$-Weyl module, which is by Theorem~\ref{sl2isobek} isomorphic to the $\lie{sl}_2\otimes \C[t]$-Demazure module $D(1,\langle\lambda,\beta^{\vee}\rangle\omega)$. Since $w_{\lambda}$ is a cyclic generator for $V$ and satisfies obviously the defining relations of $W(\langle\lambda,\beta^{\vee}\rangle\omega)$ we obtain by Lemma~\ref{isosl} a surjective homomorphism:
$$W(\langle\lambda,\beta^{\vee}\rangle\omega)\cong D(1,\langle\lambda,\beta^{\vee}\rangle\omega)\twoheadrightarrow V\subseteq W^{\Gamma}(\lambda).$$
In particular, $w_{\lambda}$ satisfies the defining relations of $D(1,\langle\lambda,\beta^{\vee}\rangle\omega)$, which contain the relation $$(x^{-}\otimes t^s)^{\max\{ 0 ,\langle\lambda, \beta^{\vee}\rangle-s \}+1}.v=0\ \forall s\in \N,$$ therefore again by Lemma~\ref{isosl} we obtain $$(X^{-}_{\beta,0} \otimes t^{ms})^{\max\{ 0 ,\langle\lambda, \beta^{\vee}\rangle-s \}+1}.w_{\lambda} = 0$$
Now suppose $\beta$ is a short root and consider the $\lie{sl}_{2,\beta}\otimes \C[t]$-submodule $V=U(\lie{sl}_{2,\beta}\otimes \C[t]).w_{\lambda}\subseteq W^{\Gamma}(\lambda).$ By the same reasons as above and Lemma~\ref{isosl} we get an surjective homomorphism $$W(\langle\lambda,\beta^{\vee}\rangle\omega)\cong D(1,\langle\lambda,\beta^{\vee}\rangle\omega)\twoheadrightarrow V\subseteq W^{\Gamma}(\lambda),$$ and therefore $w_{\lambda}$ satisfies again the relations of $D(1,\langle\lambda,\beta^{\vee}\rangle\omega)$. Using the isomorphism in Lemma~\ref{isosl} we obtain: $$(X^{-}_{\beta,j} \otimes t^{ms+j})^{\max\{ 0 ,\langle\lambda, \beta^{\vee}\rangle-(ms+j) \}+1}.w_{\lambda} = 0,\ \forall s\in\N, 0\leq j\leq m-1,$$ which proves (\ref{pot2}).\\
To prove the theorem it remains to consider the case where $\Lhg$ is of type $A^{(2)}_{2l}$. We have $(\lambda,1/2,0)\in X$, in  particulary we have $D(1/a^{\vee}_0,\lambda)=V_{\omega_0t_{\lambda-\omega_l}}(\Lambda_l).$
%Assume that $\lambda=m_1\omega_1+\cdots+m_l\omega_l$, such that $m_l-1$ is even, whereby the fundamental weights in case $B_l$ are given by $$\omega_i=\al_1+2\al_2+\cdots+(i-1)\al_{i-1}+i(\al_{i+1}+\cdots+\al_l),\ i<l, \ \omega_l=1/2(\al_1+2\al_2+\cdots+l\al_{l}).$$ By (\ref{weylgrouptrans1}) we obtain, that the affine Weyl group is given by  $W_0\semi t_{\sum_{i=1}^{l}\Z \al_i}$ and therefore we get obviously $\lambda'=m_1\omega_1+\cdots+(m_l-1)\omega_l\in \sum_{i=1}^{l}\Z \al_i$, which gives evidence of  $D(1/a^{\vee}_0,\lambda)=V_{\omega_0t_{\lambda'}}(\Lambda_l).$ 
In order to use again Corollary~\ref{demazurerelations} we reformulate (\ref{pot1}) into \begin{equation}\label{pot3}(X^{-}_{\beta,j} \otimes t^{ms+j})^{k_{\beta}+1}.w_{\lambda} = 0,\ s\ge 0,\ k_{\beta} =
\left\{\begin{array}{cl} \max\{ 0 ,\langle\lambda, \beta^{\vee}\rangle-(ms+j) \}, & \mbox{if $\beta$ is long}\\ \max\{ 0 ,\langle\lambda, \beta^{\vee}\rangle-2(ms+j) \}, & \mbox{if $\beta$ is short}\\ \max\{ 0 ,\langle\lambda, \beta^{\vee}\rangle-1/2(ms+1) \}, & \mbox{if $\beta=2\alpha$, $\al$ is short} \end{array}\right.\end{equation} 
We will prove case by case that the generator of $W^{\Gamma}(\lambda)$ satisfies the relations in (\ref{pot3}). For long roots the proof is similar to the other cases by using Lemma~\ref{isoa2}. So let $\beta$ be a short root and $\langle X_{\beta,j}^{\pm}\otimes t^{ms+j},X_{2\beta,1}^{\pm}\otimes t^{ms+1},h_{\beta,j}\otimes t^{ms+j}\rangle_{\C}$ the Lie algebra which is isomorphic to $\mathfrak{C}(A^{(2)}_{2})$ by Lemma~\ref{isoa2}. 
We consider the submodule $U(\mathfrak{C}(A^{(2)}_{2})).w_{\lambda}\subseteq W^{\Gamma}(\lambda)$, which is trivially a quotient of the $A^{(2)}_{2}$-Weyl module $W^{\Gamma}(\langle\lambda,\beta^{\vee}\rangle\omega).$ In Section~\ref{section7} Theorem~\ref{mainthma2} we prove (independent of Section 1-6) that $W^{\Gamma}(\langle\lambda,\beta^{\vee}\rangle\omega)\cong D(1/2,\langle\lambda,\beta^{\vee}\rangle\omega)$. 
The proof is finished with the observation, that the defining relations for $A^{(2)}_{2}$-Demazure module $D(1/2,\langle\lambda,\beta^{\vee}\rangle\omega)$ contain the relations $$(X^{-}_{\beta,j} \otimes t^{ms+j})^{\max\{ 0 ,\langle\lambda,\beta^{\vee}\rangle-2(ms+1) \}+1}.w=0,\ (X^{-}_{2\beta,1} \otimes t^{ms+1})^{\max\{ 0 ,1/2(\langle\lambda,\beta^{\vee}\rangle-(ms+1)) \}+1}.w=0.$$
\endproof
\end{thm}

\subsection{Fundamental Weyl modules}
In the previous section we have seen that Weyl modules are isomorphic to certain Demazure modules. Since most of the Demazure modules have a nice tensor product decomposition, see Theorem~\ref{tensorzerlegung}, we can transfer this result to most Weyl modules (only the $A^{(2)}_{2l}$ case needs more work). Using this decomposition, to compute the dimension and character of Weyl modules it is enough to describe the $\Lgg_0$ decomposition of \textit{fundamental Weyl modules} $W^{\Gamma}(\omega_i)$.
\begin{thm}\label{fundamentalzer}
Let $\omega_1,\cdots,\omega_l$ be the fundamental weights in $P^{+}_{0}$. Viewed as a $\Lgg_0$-module the fundamental Weyl modules decomposes into the direct sum of irreducible $\Lgg_0$-modules as follows:

$\bullet$ if $\Lhg$ is of type $A^{(2)}_{2l}$ 
\begin{eqnarray*}&&W^{\Gamma}(\omega_i)\cong V(\omega_i),\\
&&W^{\Gamma}(2\omega_l)\cong V(2\omega_l)\end{eqnarray*}

$\bullet$ if $\Lhg$ is of type $A^{(2)}_{2l-1}$
$$W^{\Gamma}(\omega_i)\cong \bigoplus_{s_{\bar{i}}+\cdots +s_i=1} V(s_{\bar{i}}\omega_{\bar{i}}+\cdots+s_{i-2}\omega_{i-2}+s_i\omega_i), \mbox{where $\bar{i}\in \{0,1\}$ and $i=\bar{i} \mod 2$}$$

$\bullet$ if $\Lhg$ is of type $D^{(2)}_{l+1}$
\begin{eqnarray*}&&W^{\Gamma}(\omega_i)\cong \bigoplus_{s_1+\cdots +s_i \leq 1} V(s_1\omega_1+\cdots+s_i\omega_i),\ i\neq l\\
&&W^{\Gamma}(\omega_l)\cong V(\omega_l)\end{eqnarray*}

$\bullet$ if $\Lhg$ is of type $E^{(2)}_{6}$
\begin{eqnarray*}&&W^{\Gamma}(\omega_1)\cong \bigoplus_{s\leq 1}V(s\omega_1)\\
&&W^{\Gamma}(\omega_2)\cong V(0)\oplus V(\omega_1)^{\oplus2}\oplus V(\omega_2)\oplus V(\omega_4)\\
&&W^{\Gamma}(\omega_3)\cong V(0)^{\oplus2}\oplus V(\omega_1)^{\oplus4}\oplus V(\omega_2)^{\oplus3}\oplus V(\omega_4)^{\oplus3}\oplus V(2\omega_1)\oplus V(\omega_1+\omega_4)\oplus V(\omega_3)\\
&&W^{\Gamma}(\omega_4)\cong \bigoplus_{s_1+s_4\leq 1}V(s_1\omega_1+s_4\omega_4)\end{eqnarray*}

$\bullet$ if $\Lhg$ is of type $D^{(3)}_{4}$
\begin{eqnarray*}&&W^{\Gamma}(\omega_1)\cong V(0)\oplus V(\omega_1)\oplus V(\omega_2)^{\oplus 2}\\
&&W^{\Gamma}(\omega_2)\cong \bigoplus_{s\leq 1} V(s\omega_2)\end{eqnarray*}
\proof

If $\Lhg$ is of type $A^{(2)}_{2l-1}$ or $D^{(2)}_{l+1}$ the decomposition rule is immediate from Theorem~\ref{mainthm} and Theorem 2 in \cite{FoL06}. By same reasons the theorem is true for $i=2$ if $\Lhg$ is of type $D^{(3)}_{4}$ and for $i=1,4$ in type $E^{(2)}_{6}$. For $i=1$ one can check $t_{-w_1}=w_0s_0s_2s_1s_2s_0$ and therefore with the Demazure character formula we get 
$$D_{t_{-w_1}}(e^{\Lambda_0})=D_{w_0}(e^0+2e^{\omega_2}+e^{\omega_1})\Rightarrow W^{\Gamma}(\omega_1)\cong V_{t_{-\omega_1}}(e^{\Lambda_0})\cong_{G_2}V(0)\oplus V(\omega_1)\oplus V(\omega_2)^{\oplus 2},$$ 
which proves the claim for type $D^{(3)}_{4}$.\\
So it remains to consider the nodes $i=2,3$ in type $E^{(2)}_{6}$ and the general case in type $A^{(2)}_{2l}$. In \cite{CM06}  Kirillov-Reshetikhin modules $KR(s\omega_i)$ respectively $KR^{\sigma}(s\omega_i)$ for the twisted version are defined. By inspecting the defining relations it follows that KR-modules of level $1$ (e.g. $s=1$) are precisely fundamental Weyl modules, in particular 
$$W(\omega_i)\cong KR(\omega_i) \mbox{ and } W^{\Gamma}(\omega_i)\cong KR^{\sigma}(\omega_i).$$ 
Since the decomposition of KR-modules are known as $\Lgg$ respectively $\Lgg_{0}$-modules (see \cite{KL98},\cite{C01},\cite{HKOTT02} or \cite{CM06} for instance) we obtain the predicted decomposition for $i=2,3$ in type $E^{(2)}_{6}$ and for the general case in type $A^{(2)}_{2l}$.\\ 
It remains to consider $W^{\Gamma}(2\omega_l)$, so let $\langle X^{\pm}_{\al_l,j}\otimes t^{ms+j}, X^{\pm}_{2\al_l,1}\otimes t^{ms+1}, h_{\al_l,j}\otimes t^{ms+j}|s\in \N, 0\leq j \leq m-1\rangle_{\C}$ be the Lie algebra which is by Lemma~\ref{isoa2} isomorphic to $\mathfrak{C}(A^{(2)}_{2})$. Then we obtain a surjective homomorphism $$W^{\Gamma}(2\omega)\twoheadrightarrow U(\mathfrak{C}(A^{(2)}_{2})).w_{2\omega_l}\subseteq W^{\Gamma}(2\omega_l).$$ In Section~\ref{section7} we will show that the $A^{(2)}_{2}$-Weyl module $W^{\Gamma}(2\omega)$ is an irreducible $\lie{sl}_2$-module and hence $(X^{-}_{\al_l,0}\otimes t^{2}).w_{2\omega_l}=(X^{-}_{2\al_l,1}\otimes t).w_{2\omega_l}=(X^{-}_{\al_l,1}\otimes t).w_{2\omega_l}=0$. So $W^{\Gamma}(2\omega_l)$ is isomorphic to the Kirillov-Reshetikhin module $KR^{\sigma}(2\omega_l)$, hence the decomposition is known by \cite{CM06}.
\endproof
\end{thm}
Such a similar decomposition is already known for the untwisted fundamental Weyl modules $W(\omega_i)$, see \cite{C01} or \cite{FoL06} for instance. This fact motivates us to compare the dimension of twisted and untwisted fundamental Weyl modules. For notational reasons, we have to extend certain linear functions $\lie h_0 \longrightarrow \C$ to functions on $\lie h$. So let $\mu\in P_0^+$ (with $\mu(\al_l^{\vee})\in 2\Z_{\geq0}$ if $\Lgg$ is of type $A_{2l}$). We define the extension, by abuse of notation also denoted by $\mu$, on a basis of $\lie h$ by:
$$\mu(h_i)=\left\{\begin{array}{cl} \mu(\al_i^{\vee}) & \mbox{ if $\Lgg$ is not of type $A_{2l}$}\\
0 & \mbox{ if $i\notin I_0$}\\
 (1-\frac{\delta_{i,l}}{2})\mu(\al_i^{\vee}) & \mbox{ if $\Lgg$ is of type $A_{2l}$}
\end{array}\right.$$

Since there might be a confusion in notation in the $A_{2l}$ and the $l$-th fundamental weight case only, we will use this identification in the remaining of the paper without further comment.

\begin{lem}\label{dimcheck}
Let $\omega_1,\cdots,\omega_l$ be the fundamental weights in $P_0^{+}$. We set $\epsilon=(1+\delta_{i,l})$ if $\Lgg$ is of type $A_{2l}$ and $\epsilon=1$ else, then we obtain $$\dim W^{\Gamma}(\epsilon\omega_{i})=\dim W{(\epsilon\omega_i)},\ 1\leq i \leq l.$$
\proof
Using Theorem~\ref{fundamentalzer}, Theorem 2 in \cite{FoL06} and Lecture 24 in \cite{FH91}, we obtain the following straightforward calculations:\\
$\bullet$ if $\Lhg$ is of type $A^{(2)}_{2l}, (\Lgg,\Lgg_0)=(A_{2l},B_l)$:
$$\dim W^{\Gamma}(\epsilon\omega_{i})=\binom{2l+1}{i}=\dim(V_{\Lgg}(\omega_i))=\dim W(\epsilon\omega_i)$$
$\bullet$ if $\Lhg$ is of type $A^{(2)}_{2l-1}, (\Lgg,\Lgg_0)=(A_{2l-1},C_l)$:
$$\dim W^{\Gamma}(\omega_i)=\binom{2l}{\bar{i}}+\sum_{j=1}^{\frac{i-\bar{i}}{2}}\binom{2l}{
\bar{i}+2j}-\binom{2l}{\bar{i}+2j-2}=\binom{2l}{i}=\dim W(\omega_i)$$
$\bullet$ if $\Lhg$ is of type $D^{(2)}_{l+1}, (\Lgg,\Lgg_0)=(D_{l+1},B_l)$:
\begin{eqnarray*}\dim W^{\Gamma}(\omega_i)&&=\left\{\begin{array}{cl} 2^i, & \mbox{if $i=l$}\\ 1+\sum_{j=1}^{i}\binom{2l+1}{j}, & \mbox{$i\neq l$}
\end{array}\right.=\left\{\begin{array}{cl} 2^i, & \mbox{if $i=l$}\\ \sum_{j=0}^{\frac{i-p_i}{2}}\binom{2l+2}{p_i+2j}, & \mbox{$i\neq l$}
\end{array}\right.\\&&=\left\{\begin{array}{cl} \dim V_{\Lgg}(\omega_l),\mbox{ if $i=l$}\\ \dim(V_{\Lgg}(\omega_i)\oplus V_{\Lgg}(\omega_{i-2})\oplus\cdots\oplus V_{\Lgg}(\omega_{p_i})),  \mbox{ $i\neq l$}\end{array}\right.=\dim W(\omega_i)\end{eqnarray*}
$\bullet$ if $\Lhg$ is of type $E^{(2)}_{6}, (\Lgg,\Lgg_0)=(E_6,F_4)$:
\begin{eqnarray*}&&\dim W^{\Gamma}(\omega_1)=27=\dim V_{\Lgg}(\omega_1)=\dim W(\omega_1)\\
&&\dim W^{\Gamma}(\omega_2)=378=\dim(\bigoplus_{s_2+s_6=1}V_{\Lgg}(s_2\omega_2+s_6\omega_6))=\dim W(\omega_2)\\
&&\dim W^{\Gamma}(\omega_3)=3732=\dim(V_{\Lgg}(0)\oplus V_{\Lgg}(\omega_4)^{\oplus2}\oplus V_{\Lgg}(\omega_1+\omega_6)\oplus V_{\Lgg}(\omega_3))=\dim W(\omega_3)\\
&&\dim W^{\Gamma}(\omega_4)=79=\dim (\bigoplus_{s_4\leq1}V_{\Lgg}(s_4\omega_4))=\dim W(\omega_4)\end{eqnarray*}
$\bullet$ if $\Lhg$ is of type $D^{(3)}_{4}, (\Lgg,\Lgg_0)=(D_4,G_2)$:
\begin{eqnarray*}&&\dim W^{\Gamma}(\omega_1)=29=\dim (V_{\Lgg}(\omega_1)\oplus V_{\Lgg}(0))=\dim W(\omega_1)\\
&&\dim W^{\Gamma}(\omega_2)=8=\dim (V_{\Lgg}(\omega_2))=\dim W(\omega_2)\end{eqnarray*}
\endproof
\end{lem}

%%%%%%%%%%%%%%%%%%%%%%%%%%%%%%%%%%%%%%%%%%%%%%%%%%%
%%%%%%Connection between twisted and untwisted Weyl modules
%%%%%%%%%%%%%%%%%%%%%%%%%%%%%%%%%%%%%%%%%%%%%%%%%%%

\section{Connection between twisted and untwisted Weyl modules}\label{section6}
In this section we will show that the Weyl modules $W^{\Gamma}(\lambda)$ can be realized as associated graded modules of certain untwisted Weyl modules for the loop algebra $\Lgg\otimes \C[t,t^{-1}]$. So consider for $a\in \C^{*}$ the Lie algebra homomorphism $\varphi_a$ defined as follows:
$$
\varphi_a:\Lgc\longrightarrow\Lgc,\quad x\otimes t^m\mapsto x\otimes (t+a)^m.
$$
For a $\Lgc$-module $W$ we denote by $W_a$ be the module obtained by pulling back $W$ through $\varphi_a$, i.e. $x \otimes t^s$ acts by $x \otimes (t+a)^s$. Further we denote by $\overline{W}$ be the module $W$ considered as a $\Cg$-module, obtained by the embedding 
$$\Cg\hookrightarrow \Lgg\otimes \C[t].$$
We will prove:
\begin{thm}\label{realofweyl0}
Let $\lambda=\sum^{l}_{i=1}m_i\omega_i$ be a dominant $\Lgg_0$-weight. If $\Lhg$ is a twisted Kac-Moody algebra not of type $A^{(2)}_{2l}$ we get an isomorphism of $\Cg$-modules: 
$$W^{\Gamma}(\lambda)\cong \hbox{\rm gr}(\overline{W_{a}(\lambda)}).$$
If $\Lhg$ is of type $A^{(2)}_{2l}$ and $\lambda=\lambda_1+\lambda_2\in P_0^+$, such that $m_l$ and $\lambda_2(\al_l^{\vee})$ are odd we get an isomorphism of $\Cg$-modules: 
$$W^{\Gamma}(\lambda)\cong \hbox{\rm gr}(\overline{W_{a}(\lambda_1)}\otimes W^{\Gamma}(\lambda_2)).$$
\proof
Let $\Lhg$ be not of type $A^{(2)}_{2l}$, by combining \cite{FKKS11}  and \cite{CFS08} it follows, that $\overline{W_{a}(\lambda)}$ is a cyclic $\Cg$ module. Therefore the associated graded is again cyclic and it remains to observe, that the image of the highest weight generator $\overline{\bold{w}}$ satisfies for $j\in \{0,\ldots,m-1\}$ and $h_j\in \Lh_j$ the relations
$$(h_j\otimes t^{ms+j}).\overline{\bold w}=0, \quad (s,j)\neq (0,0)$$
$$(h_0\otimes 1)\overline{\bold w}=\lambda(h_0)\overline{\bold w}.$$ 
Thus we obtain a surjective homomorphism 
\begin{equation}\label{surgh}W^{\Gamma}(\lambda)\twoheadrightarrow \hbox{\rm gr}(\overline{W_{a}(\lambda)}).\end{equation} 
In order to compare the dimension of these modules we exploit the tensor product decomposition of $W^{\Gamma}(\lambda)$ as a $\Lgg_{0}$-module by  combining Theorem~\ref{mainthm} and Proposition~\ref{tensorzerlegung}. We obtain the following :
\begin{equation}\label{tenzer1}W^{\Gamma}(\lambda)\cong W^{\Gamma}{(\omega_1)}^{\otimes m_1}\otimes \cdots \otimes W^{\Gamma}{(\omega_l)}^{\otimes m_l} \mbox{ as } \lie g_0 \mbox{-modules}.\end{equation}
An analogue decomposition was proven in \cite{FoL07} for untwisted Weyl modules for the current algebra of a simply-laced simple Lie algebra and is generalized in \cite{Na11} for the non simply-laced case. From this it follows immediately 
$$\dim \hbox{\rm gr}(\overline{W_{a}(\lambda)})=\dim W(\lambda)=\prod^{l}_{i=1}(\dim W(\omega_i))^{m_i}.$$
Hence by Lemma~\ref{dimcheck} we check that (\ref{surgh}) is in fact an isomorphism.\\

From now on, we assume that $\Lhg$ is of type $A^{(2)}_{2l}$. Since $\overline{W_{a}(\lambda_1)}$ and $W^{\Gamma}(\lambda_2)$ are cyclic $\Cg$-modules it follows with the usual arguments of \cite{FKKS11} and the Chinese remainder theorem, that the tensor product is cyclic as well. Therefore we obtain similar to (\ref{surgh}) a surjective homomorphism \begin{equation}\label{surgh2}W^{\Gamma}(\lambda)\twoheadrightarrow \hbox{\rm gr}(\overline{W_{a}(\lambda_1)}\otimes W^{\Gamma}(\lambda_2)).\end{equation}With the aim to compare the dimension on both sides we notice $$\dim \hbox{\rm gr}(\overline{W_{a}(\lambda_1)}\otimes W^{\Gamma}(\lambda_2))=\dim W(\lambda_1)\dim W^{\Gamma}(\lambda_2)=\prod^{l}_{i=1}(\dim W(\omega_i))^{\lambda_1(\al_i^{\vee})}\dim W^{\Gamma}(\lambda_2).$$ Our goal now is to prove the following tensor product decomposition: \begin{equation}\label{tenzer2}W^{\Gamma}(\lambda)\cong_{\Lgg_0} W^{\Gamma}{(\omega_1)}^{\otimes m_1}\otimes \cdots \otimes W^{\Gamma}{(\omega_{l-1})}^{\otimes m_{l-1}}\otimes W^{\Gamma}{(2\omega_l)}^{\otimes k-1}\otimes W^{\Gamma}{(\omega_l)},\end{equation}where $m_l=2k-1$
since the proposition is a immediate consequence of (\ref{tenzer2}) and Lemma~\ref{dimcheck}. To prove (\ref{tenzer2}) we investigate the character of $W^{\Gamma}(\lambda)$. By Theorem~\ref{mainthm} and Theorem~\ref{demazurcharacterformula} we obtain $$\charc W^{\Gamma}(\lambda)=\charc V_{\omega_0t_{\lambda-\omega_l}}(\Lambda_l)=D_{\omega_0t_{\lambda-\omega_l}}(e^{\Lambda_l}).$$
Suppose that $V(\mu)$ is a irreducible $B_l$-module, such that the coefficient $n_l$ is even, whereby $\mu=\sum_{i=1}^{l}n_i\omega_i$. The first step will be to show that $\charc V(\mu)$ is stable under the Demazure operators $D_i$, $i=0,\ldots,l$. The character of a finite dimensional $\Lgg_0$-module is stable under the Weyl group W and hence stable under $D_i$, $i=1,\ldots,l$. It remains to consider the case $i=0$. Note that $\al_0=\delta-2\overline{\theta}=\delta-\theta$ where $\overline{\theta}=\al_1+\cdots+\al_l$ is the highest short root of $B_l$. We define maps $s_{\overline{\theta}}:(\Lh_0)^{*}\rightarrow (\Lh_0)^{*},\ s_{\overline{\theta}}(\lambda)=\lambda-\lambda(\overline{\theta}^\vee)\overline{\theta}$ and $s_{\theta}:(\Lh_0)^{*}\rightarrow (\Lh_0)^{*},\ s_{\theta}(\lambda)=\lambda-\lambda(\theta^\vee)\theta$. Since $\overline{\theta}^\vee=2(\al_1^{\vee}+\ldots+\al_{l-1}^{\vee})+\al_l^{\vee}$ and $\theta^\vee=\al_1^{\vee}+\ldots+\al_{l-1}^{\vee}+\frac{1}{2}\al_l^{\vee}$ we get clearly $s_{\overline{\theta}}=s_{\theta}$. Thus $\nu$ is a weight in $V(\mu)$ if and only if $s_{\theta}(\nu)$ is a weight. Assume $\nu\in(\Lh_0)^{*}$ is a weight, hence $\nu=\mu-Q^{+}_0$ and therefore $\langle \nu,\al_0^{\vee}\rangle=\langle \nu,(\delta-\theta)^{\vee}\rangle=\langle \nu,-\theta^{\vee}\rangle\in\Z$. We have proved that $D_{0}$ can be defined on $\charc V(\mu)$ and $D_0=D_{-\theta}$. We obtain $$D_0(\charc V(\mu))=D_{-\theta}(\charc V(\mu))=\charc V(\mu)$$
In a second step we prove that the characters are the same by using induction on $\sum_{i=1}^{l-1}m_i+(k-1)$. So if the sum is 1 we have to show $$D_{\omega_0t_{\omega_i}}(e^{\Lambda_l})=\charc W^{\Gamma}(\omega_i+\omega_l)=e^{\frac{1}{2}\Lambda_0}\charc (V_{\Lgg_0}(\omega_i)\otimes V_{\Lgg_0}(\omega_l)),\ i<l $$$$ D_{\omega_0t_{2\omega_l}}(e^{\Lambda_l})=\charc W^{\Gamma}(3\omega_l)=e^{\frac{1}{2}\Lambda_0}\charc (V_{\Lgg_0}(2\omega_l)\otimes V_{\Lgg_0}(\omega_l)).$$
In other words, we have to figure out the $\Lgg_0$-module decomposition of $W^{\Gamma}(\omega_i+\omega_l)$ respectively $W^{\Gamma}(3\omega_l)$. By Lemma~\ref{weylprop}(2) we already know that there exists such a decomposition and since the modules are finite--dimensional every $\Lgg_0$-submodule is a direct summand. So our assignment is to find all highest weight vectors, first beginning with the highest weight vectors living in $W^{\Gamma}(\omega_i+\omega_l)[1]$. Suppose $\alpha\in\Phi_1$, such that $(X^{-}_{\alpha,1}\otimes t).w$ is a highest weight vector, i.e. the element is non-zero and the upper triangular part of $\Lgg_0$ acts by zero. We want to restrict the choice of $\alpha$ to one possible case. Note that $\alpha$ is of the form $\al_j+\cdots+\al_l$ or $2(\al_j+\cdots+\al_l)$,\ $1\leq j \leq l$ or of the form $\al_p+\cdots+\al_q$ respectively $\al_p+\cdots+\al_{q-1}+2(\al_q+\cdots+\al_l)$, $p,q \leq l-1$. If $\alpha$ is a short root, we obtain from Lemma~\ref{isoa2} \begin{equation}\label{1234}W^{\Gamma}(\langle \omega_i+\omega_l,\alpha^{\vee}\rangle\omega) \twoheadrightarrow U(\langle X_{\al,j}^{\pm}\otimes t^{ms+j},X_{2\al,1}^{\pm}\otimes t^{ms+1},h_{\al,j}\otimes t^{ms+j}\rangle_{\C}\cong\mathfrak{C}(A^{(2)}_{2})).w,\end{equation}whereby $W^{\Gamma}(\langle \omega_i+\omega_l,\alpha^{\vee}\rangle\omega)$ is the Weyl module for type $A^{(2)}_{2}$. So if $j>i$ in the representation of $\alpha$ as a sum of simple roots we get $\langle \omega_i+\omega_l,\alpha^{\vee}\rangle=1$. In Section~\ref{section7} it is shown that $W^{\Gamma}(\omega)$ is irreducible and therefore $(X^{-}_{\alpha,1}\otimes t).w=(X^{-}_{2\alpha,1}\otimes t).w=0$. Now assume $j<i$ and $(X^{-}_{\alpha_j+\cdots+\alpha_l,1}\otimes t).w\neq 0$ is a highest weight vector. Hence $0=(X^{+}_{\al_j+\cdots+\al_{i-1},0}\otimes 1)(X^{-}_{\alpha_j+\cdots+\alpha_l,1}\otimes t).w=(X^{-}_{\alpha_i+\cdots+\alpha_l,1}\otimes t).w,$ which is a contradiction to (\ref{surgh2}). In almost the same manner one sees that $(X^{-}_{2(\alpha_j+\cdots+\alpha_l),1}\otimes t).w$ cant't be a highest weigth vector. If $\alpha$ is a long root, we get with Lemma~\ref{isoa2} $$W(\langle \omega_i+\omega_l,\alpha^{\vee}\rangle\omega) \twoheadrightarrow U(\lie{sl}_{2,\alpha}\otimes\C[t]).w,$$whereby $W^{\Gamma}(\langle \omega_i+\omega_l,\alpha^{\vee}\rangle\omega)$ is the Weyl module for the current algebra $\lie{sl}_{2}\otimes\C[t]$. So if $\alpha=\al_p+\cdots+\al_q$ we obtain again $\langle \omega_i+\omega_l,\alpha^{\vee}\rangle\leq 1$ and therefore $(X^{-}_{\al,1}\otimes t).w=0$. Let $\alpha$ be of the form $\al_p+\cdots+\al_{q-1}+2(\al_q+\cdots+\al_l)$, such that $i\geq p$ and $(X^{-}_{\al,1}\otimes t).w$ is a non-zero highest weight vector. Therefore the upper triangular part acts by zero, especially \begin{flalign*}0&=(X^{+}_{\al_q+\cdots+\al_{l},0}\otimes 1)(X^{+}_{\al_p+\cdots+\al_{i-1},0}\otimes 1)(X^{-}_{\al_p+\cdots+\al_{q-1}+2(\al_q+\cdots+\al_l),1}\otimes t).w&\\&=(X^{+}_{\al_q+\cdots+\al_{l},0}\otimes 1)(X^{-}_{\al_i+\cdots+\al_{q-1}+2(\al_q+\cdots+\al_l),1}\otimes t).w=(X^{-}_{\alpha_i+\cdots+\alpha_l,1}\otimes t).w,&\end{flalign*} which is again a contradiction to (\ref{surgh2}). Hence the only possibility to get a highest weight vector of degree one is to apply $(X^{-}_{\alpha_i+\cdots+\alpha_l,1}\otimes t)$ on w. Clearly we have by Section~\ref{section7} $(X^{-}_{\alpha_i+\cdots+\alpha_l,1}\otimes t)^2.w=(X^{-}_{\alpha_i+\cdots+\alpha_l,1}\otimes t^{2s+1}).w=(X^{-}_{\alpha_i+\cdots+\alpha_l,0}\otimes t^{2s}).w=0$ for $s\geq1$, because in (\ref{1234}) we have $\langle \omega_i+\omega_l,(\alpha_i+\ldots+\alpha_l)^{\vee}\rangle=3$. Thus one can check that $(X^{-}_{\alpha_i+\cdots+\alpha_l,1}\otimes t).w$ satisfies the relations (\ref{weylrel}), (\ref{weylrel2}) in Definition~\ref{weyldef} and has weight $\omega_{i-1}+\omega_l$ with respect to $\Lh_0$.
Hence the calculations above show on the one hand that $(X^{-}_{\alpha_i+\cdots+\alpha_l,1}\otimes t).w$ is really a highest weight vector but on the other hand we get more than this, namely a surjective map $$W^{\Gamma}(\omega_{i-1}+\omega_l)\twoheadrightarrow U(\mathfrak{C}(A^{(2)}_{2l}))(X^{-}_{\alpha_i+\cdots+\alpha_l,1}\otimes t).w$$
Since $W^{\Gamma}(\omega_i+\omega_l)[1]\cong_{\Lgg_0}V_{\Lgg_0}(\omega_{i-1}+\omega_l)\cong U(\Lgg_0)(X^{-}_{\alpha_i+\cdots+\alpha_l,1}\otimes t).w$ we obtain $W^{\Gamma}(\omega_i+\omega_l)=U(\Lgg_0).w\oplus U(\mathfrak{C}(A^{(2)}_{2l}))(X^{-}_{\alpha_i+\cdots+\alpha_l,1}\otimes t).w\cong_{\Lgg_0}V_{\Lgg_0}(\omega_i+\omega_l)\oplus W^{\Gamma}(\omega_{i-1}+\omega_l)/I$, for some ideal $I$. Using (\ref{surgh2}) one can check that the ideal is zero and therefore by induction we prove our claim, because for $i=1$ we get \begin{flalign*}\charc (V_{\omega_0t_{\omega_1}}(\Lambda_l)&\cong V_{\omega_0s_0s_1\ldots s_l}(\Lambda_l))=D_{\omega_0}D_0\ldots D_l(e^{\Lambda_l})=D_{\omega_0}(e^{\Lambda_l}+e^{\Lambda_l-\al_l}+\cdots+e^{\Lambda_l-\al_l-\cdots- \al_1}+e^{\Lambda_l+\omega_1})&\\&=e^{\frac{1}{2}\Lambda_0}\charc(V_{\Lgg_0}(\omega_1+\omega_l)\oplus V_{\Lgg_0}(\omega_l))=e^{\frac{1}{2}\Lambda_0}\charc(V_{\Lgg_0}(\omega_1)\otimes V_{\Lgg_0}(\omega_l))&\end{flalign*}
Exactly the same way one can prove the existence of a surjective map $$W^{\Gamma}(\omega_{l-1}+\omega_l)\twoheadrightarrow U(\mathfrak{C}(A^{(2)}_{2l}))(X^{-}_{\alpha_l,1}\otimes t).w$$ Furthermore a more simple calculation shows $W^{\Gamma}(3\omega_l)[1]\cong_{\Lgg_0}V_{\Lgg_0}(\omega_{l-1}+\omega_l)\cong U(\Lgg_0)(X^{-}_{\alpha_l,1}\otimes t).w$. Hence $W^{\Gamma}(3\omega_l)\cong_{\Lgg_0} V_{\Lgg_0}(2\omega_l)\otimes V_{\Lgg_0}(\omega_l)$, which proves finally the initial step.
So let $\sum_{i=1}^{l-1}m_i+(k-1)>1$ and $m_i,\ i<l$ or $k-1$ such that one of them is bigger or equal to 1. Using Proposition~\ref{opaus} we get in the first case \begin{flalign*}D_{\omega_0t_{\lambda-\omega_l}}(e^{\Lambda_l})&=D_{t_{-\omega_i}}D_{\omega_0t_{\lambda-\omega_l-\omega_i}}(e^{\Lambda_l})&\\&=D_{t_{-\omega_i}}(e^{\frac{1}{2}\Lambda_0}\charc(V_{\Lgg_0}(\omega_1)^{\otimes m_1}\otimes\cdots\otimes V_{\Lgg_0}(\omega_i)^{\otimes m_{i}-1}\otimes\cdots\otimes V_{\Lgg_0}(2\omega_l)^{k-1}\otimes V_{\Lgg_0}(\omega_l)))&\\&=\charc(V_{\Lgg_0}(\omega_1)^{\otimes m_1}\otimes\cdots\otimes V_{\Lgg_0}(\omega_i)^{\otimes m_{i}-1}\otimes\cdots\otimes V_{\Lgg_0}(2\omega_l)^{k-1})D_{t_{-\omega_i}} (e^{\frac{1}{2}\Lambda_0}\charc(V_{\Lgg_0}(\omega_l)))&\\&=e^{\frac{1}{2}\Lambda_0}\charc(V_{\Lgg_0}(\omega_1)^{\otimes m_1}\otimes\cdots\otimes V_{\Lgg_0}(\omega_i)^{\otimes m_{i}}\otimes\cdots\otimes V_{\Lgg_0}(2\omega_l)^{k-1}\otimes V_{\Lgg_0}(\omega_l)).&\end{flalign*} In the second we obtain
\begin{flalign*}D_{\omega_0t_{2(k-1)\omega_l}}(e^{\Lambda_l})&=D_{-t_{2\omega_l}}D_{\omega_0t_{2(k-2)\omega_l}}(e^{\Lambda_l})=D_{-t_{2\omega_l}}(e^{\frac{1}{2}\Lambda_0}\charc(V_{\Lgg_0}(2\omega_l)^{\otimes k-2}\otimes V_{\Lgg_0}(\omega_l)))&\\&=\charc(V_{\Lgg_0}(2\omega_l)^{\otimes k-2})D_{-t_{2\omega_l}}(e^{\frac{1}{2}\Lambda_0}\charc(V_{\Lgg_0}(\omega_l)))=e^{\frac{1}{2}\Lambda_0}\charc(V_{\Lgg_0}(2\omega_l)^{\otimes k-1}\otimes V_{\Lgg_0}(\omega_l)).&\end{flalign*}
\endproof
\end{thm}
As an immediate consequence of Theorem~\ref{realofweyl0} and its proof we obtain explicit dimension formulas for Weyl modules. Such formulas for Weyl modules, as already mentioned, were previously known for untwisted current algebras (see \cite{CL06},\cite{FoL07} or \cite{Na11}).
\begin{cor}
Let $\lambda=\sum^{l}_{i=1} m_i\omega_i$ be a decomposition of a dominant weight $\lambda\in P^+_{0}$.
\begin{enumerate}
\item If $\Lhg$ is a twisted affine Kac-Moody algebra not of type $A^{(2)}_{2l}\ (l\geq 1)$, then $$\dim W^{\Gamma}(\lambda)=\prod^{l}_{i=1}(\dim W^{\Gamma}(\omega_i))^{m_i}=\prod^{l}_{i=1}(\dim W(\omega_i))^{m_i}.$$
\item If $\Lhg$ is of type $A^{(2)}_{2l}$ and $m_l=2k-1$, then $$\dim W^{\Gamma}(\lambda)=\prod^{l-1}_{i=1}(\dim W^{\Gamma}(\omega_i))^{m_i}(\dim W^{\Gamma}(2\omega_l))^{k-1}\dim W^{\Gamma}(\omega_l)=(\prod^{l-1}_{i=1}\binom{2l+1}{i}^{m_i})\binom{2l+1}{l}^{k-1}2^l.$$
\end{enumerate}
\end{cor}

\subsection{Constructions from arbitrary local Weyl modules}
In the previous section we investigate the connection between untwisted and twisted Weyl modules. We have seen that the twisted ones can be realized as associated graded modules of certain untwisted Weyl modules located in a single point. In this section we generalize this result using untwisted Weyl modules located in a finite number of points.\vskip4pt
Let $W^1,\cdots, W^k$ be finite--dimensional, graded and cyclic modules with cyclic vectors $w_1,\ldots,w_k$ for the current algebra and further let $W$ be a given cyclic graded $\Cg$-module (possibly trivial) with cyclic vector $w$.
\begin{prop}\label{cyc}
Let $a_i\in\C^{*}, 1\leq i \leq k$ be non-zero complex numbers, such that $a_i^m\neq a_j^m$ for $i\neq j$, then $\overline{W^1_{a_1}}\otimes\cdots\otimes \overline{W^k_{a_k}}\otimes W$ is a cyclic $U(\Cg)$-module, particulary we get $$\overline{W^1_{a_1}}\otimes\cdots\otimes \overline{W^k_{a_k}}\otimes W= U(\Cg).(\bold w\otimes w)$$
\proof
As $W^i$ are finite--dimensional and graded, there exists a sufficiently large $N_i$ such that $x\otimes t^s$ acts trivially for $s\geq N_i$. Thus the ideal $J_i:=\Lgg\otimes (t-a_i)^{N_i}\C[t]$ acts trivially on $W^{i}_{a_i}$. We define $\eta:\C^{*}\rightarrow \N$, $a_i\mapsto N_i$, then $\Supp \eta$ do not contain two points in the same $\Gamma$-orbit and therefore similar to the proof of Theorem~\ref{realofweyl0} we obtain that $\overline{W^1_{a_1}}\otimes\cdots\otimes \overline{W^k_{a_k}}$ is a cyclic $U(\Cg)$-module. The rest is a application of the Chinese remainder theorem.
\endproof
\end{prop}
\begin{rem}
We can consider arbitrary $\Lgg$-modules $V(\lambda_i), \lambda_i\in P^+,\ 1\leq i \leq k$ as graded and cyclic $\Lgg\otimes \C[t]$-modules, where the action is given by $$x\otimes f(t).v=f(0)x.v,\ x\in\Lgg,\ f\in \C[t].$$ Hence if $W^i=V(\lambda_i)$, it is already shown in \cite{Lau10} or in a more general setting of equivariant map algebras in \cite{NSS}, that the tensor product in Proposition~\ref{cyc} is irreducible. Moreover it is known that all finite--dimensional irreducible modules are tensor products of evaluation modules.
\end{rem}
In \cite{FKKS11} local Weyl modules for equivariant map algebras were defined and a tensor product property was proven. It was shown that if $W^i$ is an untwisted graded Weyl module, then $\overline{W^i_{a_i}}$ is an local Weyl module for $\Cg$ supported in the point $a_i$. The tensor product property gives that $\overline{W_{a_1}(\lambda_1)}\otimes\cdots\otimes\overline{W_{a_r}(\lambda_r)}$ is a local Weyl module for $\Cg$. It was shown that every local Weyl module of $\Cg$ can be obtained in this way. The following corollary, in $A^{(2)}_{2l}$ again the odd-case is considered only, shows that the dimension and $\Lgg_0$ character is independent of the support of the local Weyl module.

\begin{cor}\label{fusfirst}
Let $\lambda=\lambda_1+\cdots+\lambda_r$ be a decomposition of a dominant weight $\lambda\in P^+_{0}$ into dominant weights and let $a_1, \ldots, a_r \in \C^*$ s.t. $a_i ^m \neq a_j^m$ for $i \neq j$.
\begin{enumerate}
\item If $\Lhg$ is a twisted affine Kac-Moody algebra not of type $A^{(2)}_{2l}$, then we have an isomorphism of $\Cg$-modules: 
$$W^{\Gamma}(\lambda)\cong\hbox{\rm gr}(\overline{W_{a_1}(\lambda_1)}\otimes\cdots\otimes\overline{W_{a_r}(\lambda_r)})$$
\item If $\Lhg$ is of type $A^{(2)}_{2l}$ and $\lambda_i(\al_l^{\vee})\in2\Z_{\geq_0}$ for $1\leq i \leq r-1$ and $\lambda_{r}(\al_l^{\vee})$ is odd, then we get an isomorphism of $\Cg$-modules: 
$$W^{\Gamma}(\lambda)\cong\hbox{\rm gr}(\overline{W_{a_1}(\lambda_1)}\otimes\cdots\otimes\overline{W_{a_{r-1}}(\lambda_{r-1})}\otimes W^{\Gamma}(\lambda_{r}))$$
\end{enumerate}
\proof
By Proposition~\ref{cyc} the right hand side in (1) respectively (2) is cyclic. Hence it is easy to obtain a surjecive map of $\Cg$-modules, which is by Theorem~\ref{realofweyl0} clearly an isomorphism.
\endproof
\end{cor}
\begin{rem} As mentioned in the introduction, Weyl modules are defined in  \cite{FKKS11} in a more general way, with support in $\C$. And they are parametrized by finitely supported functions from $\C$ to $P^+$. With this corollary we have shown in all cases except the even case in $A^{(2)}_{2l}$, that the dimension and $\lie g_0$ character of a local Weyl module depends only on its $\lie g_0$ maximal weight and NOT on the support of its parametrizing function. Concluding one might be able to show that the global Weyl module is a free module for a certain algebra, which might be part of a forthcoming publication. 
\end{rem}
\begin{rem}
The same construction of an associated graded module out of finite--dimensional, graded and cyclic $\Lgg\otimes \C[t]$-modules is defined in \cite{FL99} and is called the fusion product. In the twisted case the same construction fails, since for this, one would need a pullback map like $$\sum_{j=0}^{m-1}(x_j\otimes t^{ms+j})\in \Cg \mapsto \sum_{j=0}^{m-1}(x_j\otimes (t+a)^{ms+j})\notin\Cg.$$ Therefore we have constructed in our results associated graded $\Cg$-modules out of modules coming from $\Lgg\otimes \C[t]$, which represent an analogue of fusion products.
\end{rem}
\subsection{Summary of the results}
As a conclusion we summarize our results:
Let $\lambda=m_1\omega_1+\cdots+m_l\omega_l$ be a dominant weight of $\Lgg_0$ and $\epsilon=0$ if $l$ is odd and $\epsilon=1$ else, then \\

$\bullet$ if $\Lhg$ is of type $A^{(2)}_{2}$ ($n$ is odd)
$$W^{\Gamma}(n\omega)\cong \hbox{\rm gr}(\overline{W(\omega_1)}^{\otimes(k-1)}\otimes W^{\Gamma}(\omega))\cong V_{s_1t_{(n-1)\omega}}(\Lambda_1)$$

$\bullet$ if $\Lhg$ is of type $A^{(2)}_{2l}$ ($m_l$ is odd)
$$W^{\Gamma}(\lambda)\cong \hbox{\rm gr}(\overline{W(\omega_1)}^{\otimes m_1}\otimes\cdots\otimes \overline{W(\omega_{l-1})}^{\otimes m_{l-1}}\otimes \overline{W(\omega_l)}^{\otimes(k-1)}\otimes W^{\Gamma}(\omega_l))\cong V_{\omega_{0}t_{\lambda-\omega_l}}(\Lambda_l)$$

$\bullet$ if $\Lhg$ is of type $A^{(2)}_{2l-1}$
$$W^{\Gamma}(\lambda)\cong \hbox{\rm gr}(\overline{W(\omega_1)}^{\otimes m_1}\otimes\cdots\otimes \overline{W(\omega_{l})}^{\otimes m_{l}})\cong\left\{\begin{array}{cl} V_{\omega_{0}t_{\lambda}}(\Lambda_0), & \mbox{if $m_1+3m_3+\cdots+(l-\epsilon)m_{l-\epsilon}$ is even}\\V_{\omega_{0}t_{\lambda-\omega_1}}(\Lambda_1), & \mbox{else} \end{array}\right.$$

$\bullet$ if $\Lhg$ is of type $D^{(2)}_{l+1}$
$$W^{\Gamma}(\lambda)\cong \hbox{\rm gr}(\overline{W(\omega_1)}^{\otimes m_1}\otimes\cdots\otimes \overline{W(\omega_{l})}^{\otimes m_{l}})\cong \left\{\begin{array}{cl} V_{\omega_{0}t_{\lambda}}(\Lambda_0), & \mbox{if $m_l$ is even}\\V_{\omega_{0}t_{\lambda-\omega_l}}(\Lambda_l), & \mbox{else} \end{array}\right.$$

$\bullet$ if $\Lhg$ is of type $E^{(2)}_{6}$
$$W^{\Gamma}(\lambda)\cong \hbox{\rm gr}(\overline{W(\omega_1)}^{\otimes m_1}\otimes\cdots\otimes \overline{W(\omega_{l})}^{\otimes m_{l}})\cong V_{\omega_{0}t_{\lambda}}(\Lambda_0)$$

$\bullet$ if $\Lhg$ is of type $D^{(3)}_{4}$
$$W^{\Gamma}(\lambda)\cong \hbox{\rm gr}(\overline{W(\omega_1)}^{\otimes m_1}\otimes\cdots\otimes \overline{W(\omega_{l})}^{\otimes m_{l}})\cong V_{\omega_{0}t_{\lambda}}(\Lambda_0)$$
%%%%%%%%%%%%%%%%%%%%%%%%%%%%%%%%%%%%%%%%%%%%%%%%%%%%%%%%%%%%%%%%%%%%%%%%%%%%%%%%%%%%%%%%%%%%%%%%%%%%%%%%%%%%%%%%%%%%%%%%%%%%%%%%%%%
%         
%%%%%%%%%%%%%%%%%%%%%%%%%%%%%%%%%%%%%%%%%%%%%%%%%%%%%%%%%%%%%%%%%%%%%%%%%%%%%%%%%%%%%%%%%%%%%%%%%%%%%%%%%%%%%%%%%%%%%%%%%%%%%%%%%%%
\section{Proofs for the type \texorpdfstring{$A^{(2)}_{2}$}{A^{(2)}_{2}}}\label{section7}
In this section our attention is dedicated to the twisted Kac-Moody algebra $A^{(2)}_{2}$. In the previous sections we claim that the results hold already for $A^{(2)}_{2}$, so to complete our work it misses to verify the follwing main result of this section. 
\begin{thm}\label{mainthma2}
Let $n$ be an odd integer, then the Weyl module $W^{\Gamma}(n\omega)$ is isomorphic to the Demazure module $D(1/2,n\omega)\cong V_{s_1t_{(n-1)\omega}}(\Lambda_1)$.
\end{thm}
\subsection{Properties of $\Ws$ and minimal powers}
\begin{lem}\label{gareq}
Let $I^{\sigma}$ be the left ideal in $U(\mathfrak{C}(A^{(2)}_{2}))$ generated by $\mathfrak{n}_j\otimes t^jC[t^m],(h_{\al,0}\otimes t^{2r}),(h_{\al,1}\otimes t^{2r-1}), r\geq 1, 0\leq j \leq m-1.$ Then for every $k\in\N_+$ there exists a non-zero scalar $c_k, \tilde{c_k}\in\C$ such that  
\begin{enumerate}
\item \begin{equation}\label{gareq2} (X^{+}_{\al,0}\otimes 1)^{2k-1}(X^{-}_{2\al,1}\otimes t)^{k}=\left\{\begin{array}{cl} c_k (X^{-}_{\al,1}\otimes t^k)  \mod I^{\sigma}, & \mbox{if k is odd}\\ c_k(X^{-}_{\al,0}\otimes t^k) \mod I^{\sigma}, & \mbox{if k is even} \end{array}\right. \end{equation}
\item \begin{equation} (X^{+}_{2\al,1}\otimes t)^{k-1}(X^{-}_{2\al,1}\otimes t)^{k}=\tilde{c_k}(X^{-}_{2\al,1}\otimes t^{2k-1})\mod I^{\sigma}\end{equation}
\end{enumerate}
\proof 
The first equation is a simple reformulation of Lemma 3.3 (iii) in \cite{CFS08}. We will prove the second equation by induction. For $k=1$ we get trivially $\tilde{c_k}=1$. Suppose that (2) is already true for all $p\leq k$, then \begin{flalign*}&(X^{+}_{2\al,1}\otimes t)^{k}(X^{-}_{2\al,1}\otimes t)^{k+1}=(X^{+}_{2\al,1}\otimes t)(X^{+}_{2\al,1}\otimes t)^{k-1}(X^{-}_{2\al,1}\otimes t)^{k}(X^{-}_{2\al,1}\otimes t)&\\&= \tilde{c_k}(X^{+}_{2\al,1}\otimes t) (X^{-}_{2\al,1}\otimes t^{2k-1})(X^{-}_{2\al,1}\otimes t)+(X^{+}_{2\al,1}\otimes t)\mathfrak{J}(X^{-}_{2\al,1}\otimes t),\mbox{ for some $\mathfrak{J}\in I^{\sigma}$}&\\&\equiv -\frac{1}{2}\tilde{c_k} (h_{\al,0}\otimes t^{2k})(X^{-}_{2\al,1}\otimes t)\mod I^{\sigma}&\\&\equiv 2 \tilde{c_k} (X^{-}_{2\al,1}\otimes t^{2k+1})\mod I^{\sigma}\end{flalign*}
\endproof
\end{lem}
\begin{cor}\label{cor1}
Let $n\in\N$, such that $n=2k$ if $n$ is even and $n=2k-1$ if $n$ is odd. Then we have
\begin {enumerate}
\item $(X^{-}_{2\al,1}\otimes t)^k w_n=0$
\item $\left\{\begin{array}{cl} (X^{-}_{\al,0}\otimes t^k) w_n=(X^{-}_{\al,1}\otimes t^{k+1}) w_n=0, & \mbox{if k is even}\\ (X^{-}_{\al,0}\otimes t^{k+1}) w_n=(X^{-}_{\al,1}\otimes t^{k}) w_n=0, & \mbox{if k is odd} \end{array}\right.$
\item $(X^{-}_{2\al,1}\otimes t^{2k-1}) w_n=0$

\end{enumerate}
\end{cor}
\proof
Clearly part (2) and (3) are deductions of Lemma~\ref{gareq} and part (1). Assume now $(X^{-}_{2\al,1}\otimes t)^k w_n$ is non-zero element in $W^{\Gamma}(n\omega)[k]$ of weight $-2k\omega$ if $n$ is even and $(-2k-1)\omega$ if $n$ is odd and recall that $W^{\Gamma}(n\omega)[k]$ is an integrable $\lie{sl}_2$-module, i.e. Proposition~\ref{weylprop} is applicable. That means $W^{\Gamma}(n\omega)[k]_{2k\omega}\neq 0$ respectively $W^{\Gamma}(n\omega)[k]_{(2k+
1)\omega}\neq 0$, but both are impossible, which proves part (1).
\endproof
\begin{cor}\label{fid}
For all $n\in\N$ the modules $W^{\Gamma}(n\omega)$ are finite--dimensional.
\proof
Proposition ~\ref{weylprop} implies that $\Ws_{\mu}\neq 0$ only if $\mu\in n\omega-Q^{+}_0$ and suppose that $$\Ws\cong \bigoplus_{\mu\in P^{+}_0}V(\mu)^{n_{\mu}}$$ is the decomposition of $\Ws$ into irreducible $\lie g_0$-modules. Note that the number of elements in $P^{+}_0$ with the property $\mu\in n\omega-Q^{+}_0$ is finite. The corollary follows if we prove that $\dim\Ws_{\mu}< \infty$, since this implies $n_\mu<\infty$. That the dimension can`t be infinity is a direct consequence of Corollary~\ref{cor1}.
\endproof
\end{cor}
As in the other cases we show that the Weyl modules are in connection with certain associated graded modules:
\begin{prop}\label{weylsl2fus}
Let $n\in\N$, such that $n=2k$ if $n$ is even and $n=2k-1$ if $n$ is odd. Then we get a surjective map respectively an isomorphism of $U(\mathfrak{C}(A^{(2)}_{2}))$-modules
$$\Ws \left\{\begin{array}{cl}\twoheadrightarrow \hbox{\rm gr}(\overline{W_{a_1}(\omega_1)}\otimes\cdots\otimes\overline{W_{a_{k}}(\omega_1)}), & \mbox{if n is even}\\\cong \hbox{\rm gr}(\overline{W_{a_1}(\omega_1)}\otimes\cdots\otimes\overline{W_{a_{k-1}}(\omega_1)}\otimes W^{\Gamma}(\omega)), & \mbox{if n is odd} \end{array}\right.$$ The map is given by $w_n\mapsto \underbrace{w_{\omega_1}\otimes\cdots\otimes w_{\omega_1}}_{k}$ if $n$ is even and $w_n\mapsto \underbrace{w_{\omega_1}\otimes\cdots\otimes w_{\omega_1}}_{k-1}\otimes w_{\omega}$ otherwise.
\end{prop}
\begin{rem}\label{rem1}
We will proof the isomorphism claimed in the odd case in Section~\ref{7.2} and remind that the surjectivity of the maps in Proposition~\ref{weylsl2fus} follows by weight reasons.
\end{rem}
\begin{cor}
We obtain,
$$\dim\Ws\geq\left\{\begin{array}{cl} 3^{\frac{n}{2}}, & \mbox{if n is even}\\ 3^{^{\lceil \frac{n}{2} \rceil}}2, & \mbox{if n is odd}\end{array}\right.$$
\end{cor} 
In Corollar~\ref{cor1} we proved that we can explicitly specify an integer, such that the elements with higher powers of t act by zero. In the next we will refute the question, if there exists a smaller integer with same property. To show this one can use the help of associated graded modules defined in Section~\ref{section6} and Proposition~\ref{weylsl2fus}.
\begin{lem}\label{minimum}
Let $n\in\N$ like in Corollar~\ref{cor1}. Then we have,
$$\left\{\begin{array}{cl} (X^{-}_{\al,0}\otimes t^{2r}) w_n\neq 0 ,(X^{-}_{\al,1}\otimes t^{2r+1})w_n \neq 0, & \mbox{for all $r<\frac{k}{2}$ if k is even}\\ (X^{-}_{\al,0}\otimes t^{2r}) w_n\neq 0 ,(X^{-}_{\al,1}\otimes t^{2s+1})w_n \neq 0, & \mbox{for all $r<\frac{k+1}{2}$, $s<\frac{k-1}{2}$ if k is odd} \\ (X^{-}_{2\al,1}\otimes t^{2r+1})w_n \neq 0, & \mbox{if $r<k-1$} \end{array}\right.$$
\end{lem}
Before we are in position to prove our main result of this section we will formulate another necessary proposition:
\begin{prop}\label{prop1}
Let $n\in \N$ as in Corollar~\ref{cor1}, then we have surjective homomorphisms $$W^{\Gamma}((n-2)\omega)\twoheadrightarrow \left\{\begin{array}{cl} U(\mathfrak{C}(A^{(2)}_{2}))(X^{-}_{\al,1}\otimes t^{k-1})w_n, & \mbox{if k is even}\\ U(\mathfrak{C}(A^{(2)}_{2}))(X^{-}_{\al,0}\otimes t^{k-1})w_n  & \mbox{if k is odd}\end{array}\right.$$$$W^{\Gamma}((n-4)\omega)\twoheadrightarrow U(\mathfrak{C}(A^{(2)}_{2}))(X^{-}_{2\al,1}\otimes t^{2k-3})w_n$$
\begin{proof}The proof is straightforward with Corollary~\ref{cor1}.
\end{proof}
\end{prop}
\subsection{Proof of Theorem~\ref{mainthma2}}\label{7.2}
\proof
Note that Proposition~\ref{weylsl2fus} is a direct consequence of Theorem~\ref{mainthma2} and the Demazure character formula (see Theorem~\ref{demazurcharacterformula}), since this provides us $$\dim \Ws=\dim V_{s_1 t_{(n-1)\omega}}(\Lambda_1)=\dim V_{(s_1s_0)^{\lceil \frac{n}{2} \rceil} s_1}(\Lambda_1)=3^{\lceil \frac{n}{2} \rceil}2.$$ We already know by Corollary~\ref{weyldemquot} that the Demazure module $D(1/2,n\omega)$ is a quotient of the Weyl module $\Ws$. So by Corollary~\ref{demazurerelations} it remains to show that the following relations holds: \begin{equation}\label{1}(X^{-}_{\al,0}\otimes t^{2r})^{max\{0,n-4r\}+1}w_n=0\end{equation}\begin{equation}\label{2}(X^{-}_{\al,1}\otimes t^{2r+1})^{max\{0,n-2(2r+1)\}+1}w_n=0\end{equation} \begin{equation}\label{3}(X^{-}_{2\al,1}\otimes t^{2r+1})^{max\{0,k-r-1\}+1}w_n=0.\end{equation} By Corollary~\ref{cor1} we can assume that the maximums are non-zero and further suppose that $(X^{-}_{\al,0}\otimes t^{2r})^{n-4r+1}w_n\neq 0$, hence $\Ws_{(n-2(n-4r+1))\omega}[2r(n-4r+1)]\neq 0$. By Proposition~\ref{prop1} and Proposition~\ref{weylprop} (1) we get that $$\Ws_{(n-2j)\omega}[l]=0$$ for all $l$ with $$l>\left\{\begin{array}{cl} (k-1)+\cdots+(k-j)=jk-\frac{j(j+1)}{2}, & \mbox{if $0\leq j\leq k $}\\ (k-1)+\cdots+(k-(n-j))=(n-j)k-\frac{(n-j)((n-j)+1)}{2}, & \mbox{if $k< j\leq n$} \end{array}\right.$$ Hence, $$2r(n-4r+1)\leq\left\{\begin{array}{cl} jk-\frac{j(j+1)}{2}, & \mbox{if $0\leq j\leq k $}\\(n-j)k-\frac{(n-j)((n-j)+1)}{2}, & \mbox{if $k< j\leq n$} \end{array}\right.,$$ with $j=(n-4r+1)$, which contradicts $2r<k.$ Exactly the same argumentation shows also (\ref{2}) and (\ref{3})\endproof
\begin{rem}
An inspection of the proof of Theorem~\ref{mainthma2} shows, that the condition, $n$ is odd, is not needed. Thus the relations (\ref{1}), (\ref{2}), (\ref{3}) holds also in $W^{\Gamma}(2k\omega)$, but it is easy that they are not enough. For instance in $W^{\Gamma}(6\omega)$ we have already $(X^{-}_{\al,0}\otimes t^{2})^2w_6=0$, while relation (\ref{1}) gives $(X^{-}_{\al,0}\otimes t^{2})^3w_6=0$.
\end{rem} 

%%%%%%%%%%%%%%%%%%%%%%%%%%%%%%%%%%%%%%%%%%%%%%%%%%%%%%%%%%%%%%%%%%%
% References
%%%%%%%%%%%%%%%%%%%%%%%%%%%%%%%%%%%%%%%%%%%%%%%%%%%%%%%%%%%%%%%%%%%

\bibliographystyle{alpha}
\bibliography{weylcurrent-biblist}
\end{document}